\documentclass[reqno]{amsart}

\usepackage[utf8]{inputenc}
\usepackage[english]{babel}

\usepackage{hyphenat}

\newtheorem{theorem}{Theorem}
\newtheorem{lemma}[theorem]{Lemma}
\newtheorem{proposition}[theorem]{Proposition}
\newtheorem{corollary}[theorem]{Corollary}
\newtheorem{remark}{Remark}

\RequirePackage{amsmath}   
\usepackage{amsfonts}
\usepackage{wasysym}
\usepackage{upgreek}

\usepackage[bbgreekl]{mathbbol}

\usepackage{graphicx}

\usepackage{pifont}
\usepackage{ushort}
\usepackage{yhmath}
\usepackage{calrsfs}
\usepackage{stmaryrd}

\newcommand{\half}{\mbox{\scriptsize$1\over2$}}

\newcommand{\quart}{\mbox{\scriptsize$1\over4$}}
\newcommand{\Ratio}[2]{\mbox{\small$#1\over#2$}}

\newcommand{\lamIII}{{\lambda^{\!I\!I\!I\!}}\mathstrut}
\newcommand{\lamIIIp}{{\lambda^{\!I\!I\!I\!}\,}}

\newcommand{\dlamIII}{\dot\lamBda}

\newcommand{\lgrp}{\!\big\lgroup\!}
\newcommand{\rgrp}{\!\big\rgroup\!}

\newcommand{\sgn}{\upepsilon}
\newcommand{\lamBda}{\nu}
\newcommand{\deltamu}{\delta\hspace{-0.2ex}\mu}

\newcommand*{\leftleftharpoons}{\,\raisebox{+0.07em}{\scalebox{0.75}{\;\rlap{\raisebox{+0.1em}{$\leftharpoonup$}}\raisebox{-0.1em}{$\leftharpoondown$}\;}}}

\newcommand{\Eq}{Eq.~$\!$}
\newcommand{\Eqs}{Eq.s~}

\newcommand{\Or}{\mbox{$\mathrm O$}}

\newcommand{\Pain}{Painlev\'e}
\newcommand{\PainIII}{Third Painlev\'e}
\newcommand{\PIIItran}{\hbox{{\sf P}$_{\mbox{\scriptsize III\,}}$-\,function}}
\newcommand{\PIIIeq}  {\hbox{{\sf P}$_\mathrm{III}\,$-\,equation}}

\newcommand{\url}[1]{{\tt #1}}

\usepackage{colontit}

\begin{document}

\title{
On existence and properties of roots of Third Painlev\'e transcendents
}

\author{S.I.\ Tertychniy
}

\begin{abstract}
Separate consideration of properties of roots of Third \Pain{} transcendents (\PIIItran{s})
is necessary due to 
irregularity the differential equation defining them
reveals on the subset of the phase space where its solution would vanish.
Application of the Hamiltonian formalism enables one to replace the mentioned second order differential equation
(Third \Pain{} equation)
by
two independent systems of
two nonlinear first order equations whose structures allow to name them coupled Riccati equations.
The existence of \PIIItran{s} vanishing at a given non-zero point then follows, all they being analytic thereat.
The set $\mathbb{Z}_2\times \mathbb{C}$ (or $\mathbb{Z}_2\times \mathbb{R}$) can be used for their indexing.
It proves also to be natural to use as an unknown the third order derivative rather than the original unknown itself.
After transformation of the corresponding differential equations to equivalent integral equations
the efficient algorithm of the constructing of approximate solutions to Third \Pain{} equation
in vicinity of their non-zero root
in the form of truncated power series is obtained.
An example of its application is given, its numerical validation presenting results in a graphical form is carried out.
The associated approximation applicable in vicinity of a pole of the corresponding \PIIItran{} is given as well.
The bounds from below for the distances between a pair of roots of a \PIIItran{}
and between a root and a pole 
representable in terms of elementary functions
are derived.

\end{abstract}



\maketitle

\section{Introduction}\label{sec::000}

The naming \PainIII{} transcendent is used as a unifying term
referring to solutions to representatives of
one among the six subfamilies of the family
of nonlinear second order ordinary differential equations called Painlev\'e equations \cite{GLS},\cite{Cla}.
The basics of their theory
were laid
by P.\ Painlev\'e \cite{Pa} and B.\ Gambier \cite{Ga}.
The equations they discovered and studied are distinguished by possessing the so called Painlev\'e property,
provided that their complete reducing to differential equations solved in known functions is not possible.
In turn,
the Painlev\'e property 
means that for all solutions of an equation revealing it
poles are their only
movable singularities, i.e.~singularities whose location depends on what a solution is considered.
Equivalent manifestation of the Painlev\'e property is the independence of locations of multi-valued singularities (branch points)
of the choice of a solution; in other words, they must be determined by the equation itself \cite{Cla}.

Except for several special cases
which will be named 
below,
a
\PainIII{} equation 
can be represented in the following form
\begin{equation}
                   \label{equation PIII'} 
    {\ddot{ \lambda }}
    =
   { {\dot{ \lambda }}^2 \over \lambda}
-{{\dot{ \lambda }}\over t}
-{\chi_{\infty}\,{\lambda}^2 \over t^2}
+{{\lambda}^3\over t^2}
+{\chi_0 \over { t}}
-{1\over{\lambda}}
\end{equation}
(see Ref.~\cite{DIL}, \Eq{}(4.30), {\it cf} also \Eq{}(34) in Ref.~\cite{Sl1996}).
Here $t$ is the free variable, either real or complex valued,  $\chi_0$ and $\chi_\infty$ are the arbitrary
constant parameters,
$\lambda \leftleftharpoons \lambda(t)$
is the unknown function
which is assumed to be sufficiently smooth in case of real $t$ and holomorphic if $t$ varies in a complex domain 
(such functions will be called here regular).
The overdots \raisebox{-.3ex}{`$\dot{\,\vphantom{.}\,}$'} and \raisebox{-.3ex}{`$\ddot{\,\vphantom{.}\,}$'} 
denote the first and second order derivatives
with respect to $t$, respectively.

It must be noted that
there also exists yet another
representation
of a Third Painlev\'e equation
which was actually given already in the seminal publication \cite{Pa} (see \Eq{}(7) in p.~38 therein) and reproduced in Ref.~\cite{Ga} (\Eq{}III in p.~4),
and which is
still widely used. It reads
%
\begin{equation}
                                           \label{equation PIII}  
\ddot{\bar{\mathstrut\ushort\lambda}}
=
  { \dot{\bar{\mathstrut\ushort\lambda}}^2 \over {\ushort\lambda} }
- { \dot{\bar{\mathstrut\ushort\lambda}}\over  {\ushort t}  }
+
  {{{\upalpha}}{\ushort\lambda}^2+{{\upbeta}} \over {\ushort t} }
+{{\upgamma}}{\ushort\lambda}^3
+{ {{\updelta}} \over  {\ushort\lambda} }
\end{equation}
({\it cf}.\ Ref.~\cite{GLS}, \textsection 29, \Eq{}($P_3$); Ref.~\cite{Cla}, \Eq{}(1.3)).
Here the underlining of symbols ${t}$ and $ {\lambda}  $ (implying underlining of dots which denote now derivatives with respect to $ \ushort t $)
represents the adaptation of notations
designated to emphasize
the distinction of the underlined ones from the original $t$ and $\lambda$ used in \Eq{}\eqref{equation PIII'}.
Notice also the distinction of the sets of constant parameters now comprising {\it four\/} quantities denoted by the symbols $\upalpha,\upbeta,\upgamma,\updelta$
instead of the two ones present in \Eq{}\eqref{equation PIII'}.

The interrelation of the equations \eqref{equation PIII'} and \eqref{equation PIII} can be deduced from the following statement
which is verified
by
straightforward computations.
                                    \begin{proposition}\label{prop::010}
Let the function $\lambda$ of the variable $t$ obey the equation \eqref{equation PIII'} and $\upgamma,\updelta$ be arbitrary nonzero constants. Then

\begin{equation}             \label{eq::30}
\begin{aligned}
&\mbox{the function } {\ushort \lambda}({\ushort t}) =   \upgamma^{-1/4}(-\updelta)^{1/4}\,t^{-1/2}\lambda(t)
\\
&\mbox{of the variable }{\ushort t} =  2\upgamma^{-1/4}(-\updelta)^{-1/4}\, t^{1/2}
\end{aligned}
\end{equation}
solves
\Eq{}\eqref{equation PIII} in which the rest two constant parameters are set as follows:
\begin{equation}            \label{eq::40}
\upalpha=-2\upgamma^{1/2}\chi_\infty,\;
 \upbeta=2(-\updelta)^{1/2}\chi_0.
\end{equation}
                               \end{proposition}
Conversely,
it is evident that,
given an arbitrary equation \eqref{equation PIII} with nonzero parameters $\upgamma$ and $\updelta$,
it is always possible
to find 
with the help of the above formulas
the unique constants $\chi_0$ and $\chi_\infty$
ensuring any prescribed values for the constants $\upalpha$ and $\upbeta$.
The equivalence of the equation \eqref{equation PIII'} with these constant parameters
to the proposed equation \eqref{equation PIII} then follows.

The  restriction $\upgamma\not=0\not=\updelta$
ensuring invertibility  of interrelations \eqref{eq::30}
can be regarded as a sign of genericity of a Third
Painlev\'e equation \eqref{equation PIII}.
Since our current goal is the analysis of a generic situation,
we may 
consider in what follows only the first representation \eqref{equation PIII'} of these equations,
bearing in mind 
that 
due to the above correspondence
all the relations established in its framework
can 
be conveyed to a subset of equations \eqref{equation PIII} of a full measure.

It is worth mentioning that, in the literature, the naming of \Eq{}\eqref{equation PIII'},
as against the name of \Eq{}\eqref{equation PIII},
is sometimes marked by the prime, e.g.\
Painlev\'e  III$'$ instead of Painlev\'e III, as in Ref.~\cite{DIL}, see p.~5 therein.
As it has been noted, solutions to such equations are named \PainIII{} transcendents.

Shortening references, we will usually write below `\PIIItran{s}' instead of the full term just used. 
For similar reason,
an equation of the form \eqref{equation PIII'} will also be named a \PIIIeq{}{}. 

As it  turned out,
Painlev\'e transcendents
play important role in the theory of monodromy preserving deformations
of certain linear second order differential equations which were found to belong to the family of Heun equations, see Ref.s~\cite{R}, \cite{SL}.
Generally speaking, such kind relationships take place for all subfamilies of the Painlev\'e family and the Heun family, 
see Ref~\cite{Sl2000}, Theorem in p.~749.
It can also be noted that
this matter had served an independent origin of interest to Painlev\'e equations
playing role of another starting point for their finding and investigation.

In our case 
the corresponding linear
equations, named  double confluent Heun equations (hereinafter referred to, for the sake of brevity, as {\sf DCHE}),
constitute a subfamily of the family of Hein equations, see Ref.~\cite{SW}.
The basic relations of \PIIIeq{}{s} and {\sf DCHE\,}s in the context of the theory of isomonodromic deformations
of the latter 
are given in Ref.~\cite{DIL}, see section 4.8 therein.
Their  substantiation realized in a straightforward way is discussed in Ref.~\cite{T}.
In particular, the principal role of zeros of \PIIItran{s} is there emphasized.

It is also important mention  a
bridge from the aforementioned pure mathematical matter
to physics, namely, to the modeling of overdamped Josephson junctions (hereinafter, Jj).
{\sf DCHE\,}s are known as the efficient tool 
used in particular for characterization of the so-called phase-lock areas in the parameter space attributed to Jj.
With a view to development of this direction,
the role  of \PIIItran{s} 
including 
their 
application in realization of
isomonodromic deformations of {\sf DCHE}s
is 
discussed
in Ref.s~\cite{BG}, \cite{Gl1}, \cite{Gl2}.

More precisely, 
these are
poles of \PIIItran{s} which
appear there a matter of primary interest. 
However, 
roots and poles of \PIIItran{s} 
are closely related in a simple way and all properties of the latter can be characterized in terms of properties of the former and vice versa.
Indeed, \PIIIeq{s}
possess some symmetries
of which we get here the one represented by the transformation $T_2$ defined by the second equation from the triplet (29.18)
given in Ref.~\cite{GLS}.
It states that 
the replacement
\begin{equation*}                              
{\ushort\lambda}(\ushort t) \! \leftleftharpoons {\ushort\lambda}(\ushort t)^{-1}
\end{equation*}
of the unknown function $ {\ushort\lambda} $ preserves the fulfillment of \Eq{}\eqref{equation PIII},
provided the constant parameters are concurrently undergone the involutive pairwise replacements
\begin{equation*}                              
\upalpha \leftleftharpoons -\upbeta,
\upbeta \leftleftharpoons -\upalpha,\;
\upgamma \leftleftharpoons -\updelta,
\updelta \leftleftharpoons -\upgamma.
\end{equation*}
Such an invariance
can be verified by straightforward computation.

A cognate property of the equation \eqref{equation PIII'}
also verifiable by straightforward computation
reads:
the replacement
\begin{equation} 
                                  \label{root2pole and vv} 
{\lambda}( t)\! \leftleftharpoons t \,\lambda(t)^{-1}
\end{equation}
of the unknown function $\lambda$ retains its fulfillment in case of concurrent involutive interchange
\begin{equation}
                                   \label{eq::60}
\chi_0 \rightleftharpoons \chi_\infty
\end{equation}
 of the values of the constant parameters $\chi_0$  and $\chi_\infty$.

The above interrelations clearly suggest
that any pole of a \PIIItran{} situated not at zero
occupies the point of
a root of a definite \PIIItran{}
characterized by
usually
distinct but known constant parameters.

\bigskip

The goal of the present work may be expressed as the consideration of some basic properties
of roots of a Third Painlev\`e transcendent beginning with issue of their existence and regularity.
Indeed, having defined these functions as solutions to  equation \eqref{equation PIII'},
one has to take into account that the latter 
is not always well defined in itself. %
In particular, in case of generic regular $\lambda(t)$, the equation right hand side is not regular at the point $t=0$.
As a rule, it also diverges at non-zero roots of unknown $\lambda(t)$ that leads to necessity of separate consideration of
behavior of \PIIItran{s} near such points.

It turns out however
that there exists
some transformation of the original \PIIIeq{} \eqref{equation PIII'}, arising as the result of application
of the Hamiltonian formalism which produces its representation in terms of the two systems of pairs of first order nonlinear
differential equations (Hamilton equations, see \Eqs\eqref{Ham eq first}, \eqref{Ham eq second} below)
such that
their
right hand sides are polynomials of the second order with respect to a one unknown function and
a linear function of the another one.
Accordingly, the classical theory of ordinary differential equations may be here applied
from which the existence of solutions regular at non-zero roots of $\lambda(t)$ directly follows.
Its subsequent consequence is the obvious feasibility
of finding of regular \PIIItran{s} vanishing at an arbitrary pre-specified point distinct from zero.

The set of all such functions is identified with the set of solutions to the
Cauchy problem for the equation \eqref{equation PIII'}
requiring of $\lambda(t)$ to take zero value at the noted point.
However, there is a difference of its treatment as compared to the standard one
used the case of  regular second order ordinary differential equations.

The systems of differential equations with solutions directly linked to analytic \PIIItran{s} vanishing
at the given point are obtained from the noted Hamilton equations by means of a change of the unknown function
transferring its role to the own third order derivative,
see \Eqs\eqref{ODE for mu with lamIII}, \eqref{ODE for lamIII} below.
The transformed equations prove to be singular at the root of \PIIItran{s} in question but still admit solutions regular thereat.
The latter arise in case of imposing certain constraint on the corresponding initial data.
A detailed consideration
leads to the following
statement (see section \ref{subs::50.40}).
                  \begin{theorem}           \label{teo::020}
Let an arbitrary $t_0\not=0$ be chosen.
Then
there exist two families of solutions $\lambda(t)$ to  \Eq{}\eqref{equation PIII'} 
each of which vanishes at the point $t_0$ and is analytic thereat.
The noted families are distinguished by the values of the derivative $\dot\lambda(t_0)$
equal to $+1$ for one of them and to $-1$ for another.
Within a single family,
a solution
is identified
by
the value of its  third derivative $\dddot\lambda(t_0)$ which may be arbitrary.
On the contrary, the values of the second derivative $\ddot\lambda(t_0)$ 
are the same for all the members of a family,
being equal to 
$\big(\dot\lambda(t_0)-\chi_0\big)/t_0$. 
                   \end{theorem}
This means, in particular, that for arbitrary chosen constant parameter values
and any  point distinct from zero
there exists
a \PIIItran{} 
vanishing at the latter and regular (analytic or real analytic) in its neighborhood.

The systems of differential equations producing solutions equivalent to \PIIItran{s} vanishing at the given point
can be converted 
to certain systems of non-linear integral equations devoid of irregularities  in domain of regularity of unknown functions
and  automatically implementing the required constraint to the initial data for the former.

The integral equations thus obtained can be used to develop an
iterative algorithm for generating their formal solution in the form of  power series in deviation of argument $t$
from the point $t_0$ where the associated \PIIItran{} is claimed to be equal to zero.
The bounds from above for  its coefficients
are derived 
with the help of the Cauchy-Bunyakovsky-Schwarz inequality
and are proven by means of the mathematical
induction. They result in the establishing of convergence of the series in question and provide an estimate from below
for the radius of the corresponding convergence area.
Thus an analytic function here appears providing, by construction, a solution to the given
\PIIIeq{} which vanishes at the point designated.
A 
summary of these conclusions 
is given in theorem \ref{teo::120} from section \ref{subs::50.40}. 

There are actually many \PIIItran{s} distinguished by the same properties as the one thus obtained.
Their set is specified as follows.
First, only the values $\dot\lambda=1$ and $\dot\lambda=-1$ the first derivative can be equal to at the root
are allowed and they both occur. 
Second, the second order derivative evaluated at a root is linked with the first order one but for the rest it is fixed.
Finally, it is allowed for the third derivative to take on any value thereat.
The latter can be used for identification of particular solutions vanishing at a given point,
although in conjunction with value (in fact, the sign) of the first derivative.
The set $\mathbb{Z}_2\times\mathbb{C}$ (or $\mathbb{Z}_2\times\mathbb{R}$)
can thus be used
for the 
1-to-1 indexing of \PIIItran{s} vanishing at the given point.

The iterative algorithm based on the integral equations
yielding solutions to the \PIIIeq{}
in case of their vanishing at a given point
can be used for generating their approximations represented by truncated power series
coinciding with leading parts of the corresponding exact infinite series. 
Computations may start from `the void approximation' constituted of identically zero functions.

In that way, an approximate representation  accurate up to the $8^{\mbox{\scriptsize th}}$ order inclusively
of a
generic \PainIII{} transcendent
vanishing at a given point
was derived, see section \ref{sec::060} below and, in particular, \Eqs\eqref{lam_6 definition}, \eqref{PIII approx}.
Therein, the numerical validation,
applying 
visualization of the relevant functions and their relationships, is carried out.
In particular, it is demonstrated that \PIIItran{} can be evaluated on the segment linking a pair of its nearby roots
by means of a proper application of two instances of the aforementioned approximation `anchored' at the latter.

Additionally, using the known relation of roots and poles of the properly related \PIIItran{s} noted above,
the approximate representation of a generic \PIIItran{} possessing not root but a pole
at the given location is displayed in section \ref{sec::070}.

The estimate of minimum for radius of convergence for series determining solution to
 \Eqs\eqref{ODE for mu with lamIII}, \eqref{ODE for lamIII} is used for constructing of the explicit bound from below
for the distance between a pair of distinct roots of a \PIIItran{}, as well as for the mutual distance separating a root and a pole.
See theorem \ref{teo::140}, corollary \ref{coro::130}.

\section{Third Painlev\'e equation from viewpoint of the Hamiltonian formalism}\label{sec::010}

It has been found that all Painlev\'e equations can be interpreted as evolution equations for certain dynamical systems described
by non-autonomous Hamiltonians
expressed in terms of polynomial functions (being also rational functions with respect to dependence on the time variable).
In case of \Eq{}\eqref{equation PIII}, the corresponding Hamiltonian was first presented in Ref.~\cite{O}.
For \Eq{}\eqref{equation PIII'}, the associated Hamiltonian can be found in Ref.~\cite{DIL}, see Eq.~(4.29).
The equivalent expression  was given in
the earlier publication \cite{Sl1996}, see \Eq{}(33) therein,
although it needs some adaptation to notations here used.

We adopt here a slightly generalized version of the Hamiltonian
for Third \Pain{} equation
presented in Ref.~\cite{DIL} defining it as follows. 
\begin{equation}
                                 \label{Hamiltonian} 
{\mathcal H}=
{\mathcal H}(\lambda,\mu, t)
= t^{-1}\big(\lambda^2\mu^2-(\lambda^2-\lambda + \sgn \,(\chi_0\lambda - t) )\mu + \half(\chi_\infty+\sgn \chi_0-1)\lambda\big).
\end{equation}
In accordance with 
terminology used in the Hamiltonian theory,
the variable $\lambda$ is named 
the generalized coordinate,
the variable $\mu$ is the generalized
momentum conjugated to $\lambda$,
and $t$ plays role of the evolution time.
As compared to the expression given in Ref.~\cite{DIL},
its modification
 leading to \eqref{Hamiltonian}
comes down to addition of yet another constant parameter $\sgn$
(the remaining symbols $\chi_\infty$ and $\chi_0$ refer to eponymous constant parameters found in  \Eq{}\eqref{equation PIII'}). 
We may regard this as a generalization since for $\sgn = +1$ the original definition referred to above results.
Let us inspect what happens in case of other choices of $\sgn$.

To that end, let us consider the Hamilton's equations
\begin{equation*}
\dot\lambda(t)=
{{
\raisebox{-1.4ex}[1ex]{\Big\lfloor}%
_{\hspace{-2.1ex}\raisebox{+0.3ex}{
 \mbox{\scriptsize
 \raisebox{0.7ex}[1ex]{
 \rlap{
 $
\lambda\leftleftharpoons\lambda(t),\; \mu\leftleftharpoons\mu(t)
 $
 } }}
}}
}
\hspace{-0.6ex}
\partial {\mathcal H}/\partial \mu
}
,\;\hspace{3em}
\dot\mu(t)=
-
{{
\raisebox{-1.4ex}[1ex]{\Big\lfloor}%
_{ \hspace{-1.8ex}   \raisebox{+0.3ex}{
 \mbox{\scriptsize
 \raisebox{0.7ex}[1ex]{
 \rlap{
 $
\lambda\leftleftharpoons\lambda(t),\; \mu\leftleftharpoons\mu(t)
 $
 } }}
 }}
}
\hspace{-1.0ex}
\partial {\mathcal H}/\partial \lambda
}.
\end{equation*}
which in the case of Hamiltonian \eqref{Hamiltonian}, omitting indication of the function argument, look as follows:
\begin{eqnarray}
                     \label{Ham eq first}   
\dot\lambda &=& t^{-1}\left(\sgn \,t - (\sgn\,\chi_0 - 1)\, \lambda + (2 \mu - 1)\, \lambda^2\right),
\\
                     \label{Ham eq second}   
\dot\mu &=& t^{-1}\big(
-\half(\chi_\infty + \sgn\,\chi_0 - 1)
+
 (\sgn\,\chi_0 - 1 + 2\lambda)\,\mu - 2\lambda\,\mu^2)
\big).
\end{eqnarray}
Here the right hand side of \Eq\eqref{Ham eq first} is a quadratic polynomial in $\lambda$ and,
similarly,
the right hand side of \Eq\eqref{Ham eq second} is a quadratic polynomial in $\mu$. Besides, they are linear in the `complementary' unknowns
$\mu$ and $\lambda$, respectively.
These features might serve justification for the naming the system \eqref{Ham eq first},\eqref{Ham eq second} the coupled Riccati equations.

Let us assume that $\lambda(t_0)=0$ for some $t_0\not=0$. Then
\Eq\eqref{Ham eq first} reduces to the equality
\begin{equation}    \label{eq::100}
\dot\lambda(t_0)=\sgn.
\end{equation}
In turn, for such $\lambda$, \Eq\eqref{equation PIII'} yields
$ 
\dot\lambda(t_0)^2=1.
$ 
Accordingly, 
the constraint 
\begin{equation}
                \label{this is sign switch} 
\sgn^2=1
\end{equation} 
arises here in a natural way. 
Its role is explained 
by the following statement 
establishing also the concordance of  Hamiltonian \eqref{Hamiltonian} with \Eq\eqref{equation PIII'}.
%
\begin{proposition}                     \label{prop::030}
The system of Hamilton's equations \eqref{Ham eq first} and \eqref{Ham eq second}
yields solutions
to \Eq\eqref{equation PIII'}, provided \Eq\eqref{this is sign switch}
is met.
\end{proposition}
There exist therefore two admissible values of $\sgn$ that is equivalent to
admissibility of only two values $+1$ and $-1$ of $\dot\lambda$ at any nonzero root of $\lambda$
where the derivatives $\dot\lambda$ and  $\ddot\lambda$ exist.
While $\sgn=+1$ leads to the well known result,
there is another opportunity $\sgn=-1$ 
yielding a formally distinct Hamiltonian which leads however to the same equation \eqref{equation PIII'}. 

It is worth mentioning that in case of consideration of solutions to \Eq\eqref{equation PIII'}
vanishing at a given point, the transition to the system of  equations \eqref{Ham eq first}, \eqref{Ham eq second}
should not be considered as an embodiment of an equivalence relation since it includes the additional discrete parameter
$\sgn$.
In framework of the problem in question,
it is the pair of such systems implementing jointly the both admissible cases $\sgn=1$ and $\sgn=-1$
which may be regarded as the equivalent to \Eq\eqref{equation PIII'}. 

Considering now the system of  equations \eqref{Ham eq first}, \eqref{Ham eq second}
with $\sgn$ obeying the condition \eqref{this is sign switch}, one can construct its solution
assigning, in advance, any values to the initial data $\{\lambda(t_0), \mu(t_0)\}$
for the Cauchy problem
(assuming however that $t_0\not=0$), and the result will consist of
analytic functions. In particular, one may also set $\lambda(t_0)=0$.
This would lead to solution analytic at $t_0$ in which the function $\lambda(t)$ is just a \PIIItran{} vanishing at $t_0$.
It
depends on the values of $\sgn$ and $\mu(t_0)$
to be
chosen in advance
and is thus not unique.

The following statement
holds therefore true albeit this does not directly follow from the original structure of the equation in question.
\begin{proposition}\label{prop::040}
Given any point $t_0\not=0$ of the complex plane,
there exist analytic solutions to  \Eq\eqref{equation PIII'} (\PainIII{} transcendents) vanishing thereat.
\end{proposition}

It can be supplemented by the following 
\begin{proposition}\label{prop::050}
Any \PainIII{} transcendent is analytic at each its non-zero root.
\end{proposition}
Indeed, in accordance with \Pain{} property which \PIIItran{s} possess any their irregularity
observed at a point distinct from zero is a simple pole.
At the same time
any non-zero root of the given \PIIItran{} gives rise to irregularity of certain \PIIItran{} observed for the same value of the argument.
Such a \PIIItran{} arises as the result of   the transformation \eqref{eq::60} interchanging the associated constant parameters
and the solution of the new \Pain{} equation, now
with pole, is obtained by means of  the transformation \eqref{root2pole and vv}.
Its form makes it evident that the pole may here result only if the original\ \PIIItran{} is no mere zero but also analytic thereat.

\section{Some features of nonzero roots of \PainIII{} transcendents}\label{sec::020}

At a nonzero root of a \PIIItran{} $\lambda$,
the second Hamilton's equation  \eqref{Ham eq first}
reduces to the constraint \eqref{eq::100}.
The first one, \Eq\eqref{Ham eq second},
takes thereat the form
$$                             
\dot\mu(t_0)
=t_0^{-1}\big( (\sgn\chi_0-1)\mu(t_0) - \half(\chi_\infty +\sgn\chi_0-1)\big).
$$
It determines $\dot\mu(t_0)$ but does not restrict $ \mu(t_0) $ to which an arbitrary value is allowed to be assigned.
Thus for both admissible choices $\sgn=1$ and $\sgn=-1$ 
the corresponding instances of
Eq{}s~\eqref{Ham eq first}, \eqref{Ham eq second}
define two families of analytic solutions
verifying also \Eq\eqref{equation PIII'}
and vanishing at $t_0$.
Their elements can be
indexed 
by a single constant, for instance, by value of $\mu(t_0)$. 

It is obvious 
that
the knowledge of the pair $\{\sgn,\mu(t_0)\}  $ is sufficient for determination  of 
which system of \Eq{}s~\eqref{Ham eq first},\eqref{Ham eq second}, dependent on the choice of $\sgn$, is used
and, then, what a solution of the latter with the second unknown $\lambda$ (\PIIItran, in fact) vanishing at $t_0$ is considered.
It is however reasonable to have at our disposal a characteristic of 
such a \PIIItran{} bearing the same data but avoiding referring to  $\mu$.
This would obviously be more natural when dealing with just \Eq\eqref{equation PIII'}
without its conversion to \Eq{}s~\eqref{Ham eq first},\eqref{Ham eq second}.

To that end, let us consider some low order terms in the Taylor expansion of the function $\lambda$ obeying \Eq\eqref{equation PIII'}.
\Eq\eqref{eq::100} implies that in the lowest order approximation $\lambda(t) = \sgn\cdot(t-t_0)+O\big((t-t_0)^2\big)$.
Having added one more term to it,
we find upon substitution of
$ 
\lambda(t) = \sgn\cdot(t-t_0)+   \half\ddot\lambda(t_0)(t-t_0)^2  +O\big((t-t_0)^3\big)
$ 
 into \Eq\eqref{equation PIII'}
that
\begin{equation*}                       
\ddot\lambda(t_0)=t_0^{-1}(\sgn - \chi_0).
\end{equation*}
The first and the second  order derivatives of a solution
turn out therefore to be fixed by the very equation
and 
hence
cannot be used for its own characterization.

Such a peculiarity
is however settled already on
the next order of the expansion accuracy.
Indeed, the representation
\begin{equation}                      \label{eq::120}
\lambda(t) = \sgn\cdot(t-t_0) +  {(\sgn - \chi_0) \over 2 t_0 }(t-t_0)^2 + (t-t_0)^3 \lamBda(t),
\end{equation}
replacing $\lambda(t)$ by another unknown function $\lamBda(t)$ also analytic at $t_0$,
yields upon substitution into \Eq\eqref{Ham eq first} the equality
\begin{equation}           \label{smoothness condition}   
2\mu(t_0)= 1 + \sgn\cdot(1-\chi_0^2)/(2t_0) + 3 t_0 \lamBda(t_0).
\end{equation}
It puts the set of quantities $\mu(t_0)$, which play role the initial data for \Eq{}s~\eqref{Ham eq first},\eqref{Ham eq second}
in case of the vanishing of $\lambda(t_0)$, into the 1-to-1 correspondence to the set of quantities $\lamBda(t_0)$
coinciding, up to the factor of $6$, with the third order derivatives $\dddot\lambda(t_0)$.
Accordingly, the initial data for solutions to  \PIIIeq{}   \Eq\eqref{equation PIII'} vanishing
at the point $t_0\not=0$ of its specification is determined by the pairs $\big\{\sgn=\pm 1,  \dddot\lambda(t_0) \big\}$
in which the second quantity may be arbitrary.

\section{Integral equations controlling 
Third Painlev\'e transcendents near their roots}\label{sec::030}

The expansion \eqref{eq::120} can be regarded as the invertible replacing of the unknown $\lambda$
by the associated function $\lamBda$.
The latter, when considered in conjunction with the subsidiary unknown $\mu$,
obeys the own equations which are an alternative albeit equivalent representation of Eq{}s~\eqref{Ham eq first}, \eqref{Ham eq second}.
The underlaid relationships
substantiated by straightforward computations
are as follows.
                  \begin{proposition}   \label{pro:020}
If a solution $\lambda(t)$
to \Eq{}\eqref{equation PIII'} is sufficiently smooth (in particular, analytic) at the point $t=t_0\not=0$
where   $\lambda(t_0)=0$
then
it yields through \Eq\eqref{eq::120}
a solution $\lamBda(t)$ to
the following system of the first order ordinary differential equations
%
\begin{eqnarray}
                                \label{ODE for mu with lamIII}    
t\,\dot\mu(t) &=& W_{\!\mu}[\mu,\lamBda](t - t_0, t),
\\
                                \label{ODE for lamIII}            
t\,\dot\lamBda(t) &=& W_{\lamBda}[\mu,\lamBda](t - t_0, t),
%
\\ &&  \hspace{-7.03em}
                                  \label{W_mu definition}            
\begin{aligned}[t] \mbox{where }
W_{\!\mu}[\mu,\lamBda](\delta t, t)
=   &\;
-\half (\chi_\infty + \sgn\,\chi_0 - 1) - (1 - \sgn\,\chi_0)\mu(t)
\\ &\;
-2\delta t\,(\mu(t)-1)
\,\mu(t)\,
\big(
\sgn + \delta t\,(\sgn - \chi_0)/(2t_0)+\delta t^2\lamBda(t)
\big),
\end{aligned}
\\&&  \hspace{-7.03em}
                                \label{W_lam definition}            
\begin{aligned}[t] \hphantom{\mbox{where }}
W_{\lamBda}[\mu,\lamBda](\delta t, t)
=   &\;
\phantom{-}\,\,
\delta{}t^{-1}\big(\sgn\,(\chi_0^2-1)/(2t_0) - 1 + 2\mu(t) - 3t_0\,\lamBda(t)\big)
\\&\;
+  t_0^{-1}(1 -\sgn\,\chi_0)(2\mu(t)-1)
 - (2 + \sgn\,\chi_0)\lamBda(t)
\\& \hspace{0.3em}
+ \delta{}t\,(2\mu(t)-1)
\left(2\sgn\,\lamBda(t)
+\left(   \left(\sgn -\chi_0\right)\!/(2t_0)  +  \delta{}t\, \lamBda(t)\right)^2 \right),
\end{aligned}
\\&& \hspace{-6.3em}
\mbox{and where }\sgn=\lambda'(t_0)\mbox{ is equal to either $+1$ or $-1$.}
\nonumber
\end{eqnarray}
Conversely, any solution to \Eqs\eqref{ODE for mu with lamIII}, \eqref{ODE for lamIII}
with any of $\sgn=\pm1$
gives rise to solution to \Eq{}\eqref{equation PIII'} vanishing at $t=t_0$
                             \end{proposition}
\hangindent=1.0em
\hspace{0em}
{\sl A note in regard to notations}: In the above formulas, as well as in similar contexts below, the list of arguments of a function
enclosed in square brackets contains the symbols of `inferior' functions involved in its unfolding 
while its `ordinary' arguments
which assume ultimately numerical values
 are shown enclosed by parentheses, as usual.
The dependence on the constant parameters $\chi_0,\chi_\infty,\sgn,t_0$ is not displayed for the sake of brevity.

\medskip

\noindent
It has to be mentioned that
in general case, i.e.~for generic values of  $\mu(t_0)$ and $\lamBda(t_0)$,
the right hand side of \Eq{}\eqref{ODE for lamIII} is not defined
at the point $t=t_0$, being also unbounded in its vicinity.
Indeed, the first line in the definition \eqref{W_lam definition}
contains the fraction with denominator tending to zero as $t\to t_0 $.
As a consequence, $\dlamIII(t)$ is singular at $t_0$, provided \Eq{}\eqref{ODE for lamIII} is fulfilled.

There exists however a subset of the set of admissible initial data
for solutions to
\Eqs{}\eqref{ODE for mu with lamIII}, \eqref{ODE for lamIII}
such that
the noted fraction
\begin{equation}
                                           \label{xi fraction}
\xi(t)= (t-t_0)^{-1}\big(\sgn\,(\chi_0^2-1)/(2t_0) - 1 + 2\mu(t) - 3t_0\,\lamBda(t)\big)             
\end{equation}
is well defined and smooth even at $t=t_0$.
They are singled out by the imposing the constraint \eqref{smoothness condition} which forces the fraction numerator to also vanish at zero of the denominator.
Assuming \Eq{}\eqref{smoothness condition} to be fulfilled,
it is seen that if $\mu,\lamBda \in C^N $ for $ N>1$ then the suspicious fraction \eqref{xi fraction} belongs to the same class of smoothness
everywhere except at $t=t_0$ where a priori only the smoothness of order $N-1$ is guarantied.
At the same time,
in case of analytic $\mu$ and $\lamBda$, \Eq\eqref{smoothness condition} makes the fraction $\xi(t)$  also analytic,
as well as the right hand side of \Eq\eqref{ODE for lamIII}  (while the right hand side of \Eq\eqref{ODE for mu with lamIII}
is  analytic for arbitrary analytic $\mu$ and $\lamBda$).
Then no incompatibility in degrees of smoothness of the right and left hand sides of
\Eqs\eqref{ODE for mu with lamIII}, \eqref{ODE for lamIII} may arise. It should also be noted that below the possibility of finite differentiability of
$\mu$ and $\lamBda$ will be excluded, so all their derivatives at $t_0$ must exist.

The constructing of a \PIIItran{} in vicinity of its root
on the base of \Eqs\eqref{ODE for mu with lamIII}, \eqref{ODE for lamIII}
requires therefore the imposing of the additional condition expressed by \Eq\eqref{smoothness condition}
besides specification of the natural initial data for the unknown functions
obeying first order differential equations,  i.e.\ the picking and fixation of values of $\mu(t_0)$ and $\lamBda(t_0)$.
Indeed, if \Eq\eqref{smoothness condition} is not fulfilled
then
solutions to \Eqs\eqref{ODE for mu with lamIII}, \eqref{ODE for lamIII}
equivalent to \Eqs\eqref{Ham eq first}, \eqref{Ham eq second}
still exist
but for them
the associated
function
$\lambda(t)$  is either nonzero or singular at $t=t_0$; in both cases $\lamBda(t)$ is obviously singular thereat as well.

The insufficiency of the standard setting of the Cauchy problem
indicated above
can be remedied 
by means of transition from differential equations
to the associated system of integral equations. In our case the latter look as follows.
\begin{eqnarray}
%
%
                                      \label{int eq for mu}            
 &&
\begin{aligned}[t]
\mu(t)=&\;
\half
\big(
    1 - \sgn\,(\chi_0^2 - 1)/(2t_0) +3\,t_0\, \lamBda(t_0) 
\big)
\\&\;
+
  (t - t_0)\,t_0^{-1}
  \big(
       - \half(\chi_\infty+\sgn\,\chi_0-1) - \mu(t)
\\&\; \hphantom{+ (t - t_0)\,t_0^{-1}  \big(   } \;
  +
   \int^1_0 \!\!
  d \sigma\:
  \Omega_{\mu}[\mu,\lamBda](\sigma,t - t_0) \big).
\end{aligned}
\\
&&
                     \label{int eq for lamIII}                         
\begin{aligned}[t]
\lamBda(t)
=-(t-t_0)t_0^{-1}\,\lamBda(t)
+
t_0^{-1}\!
\int_0^1 \!\! d \sigma\,
\sigma^2\,
\Omega_{\lamBda}[\mu,\lamBda](\sigma,t - t_0),
\end{aligned}
\end{eqnarray}
The kernels  $ \Omega_{\mu}, \Omega_{\lamBda}$ used therein 
are defined as follows.
 \begin{eqnarray}                  \label{Omega_mu definition}              
 &&\hspace{-2.5em}
 \begin{aligned}[b]\hspace{-0.1em}
\Omega_{\mu}[\mu,\lamBda](\sigma,\delta t)
=&\;
\mu(\tau)
\Big(
\sgn\,\chi_0
+ 2\,\delta t\,\sigma
\big(1 - \mu(\tau)\big)\times
\\[-2ex] &\hspace{8em}
\big(\sgn - \delta t\,\sigma\,(\chi_0-\sgn) + \delta t^2\sigma^2\lamBda(\tau)\big)
\Big),
\end{aligned}
\\
                                                        \label{Omega_lam definition}              
 &&\hphantom{.}\hspace{-2.8em}
 \begin{aligned}[b]
\Omega_{\lamBda}[\mu,\lamBda](\sigma,\delta t)
=&\;
 \sgn\, (\chi_0^2 - 1)/(2t_0) - 1 + 2\,\mu( \tau)
 \\&\;
 -\delta t\,\sgn\,\sigma\,
 \left( (\chi_0-2\sgn)\,\lamBda(\tau) +t_0^{-1}(\chi_0 - \sgn)(2\mu( \tau) - 1)\right)
 \\&\;
 +\delta t^2 \sigma^2
 (2\mu( \tau) - 1)\times
 \\
  &\; \hphantom{   +\delta t^2 \sigma^2     }\;
 \left(2\sgn\,\lamBda(\tau) + \big(  (2t_0)^{-1}(\chi_0 - \sgn) - \delta t\,\sigma\,\lamBda(\tau) \big)^2 \right),
\\[-3.2ex]
   \phantom{.}
 \end{aligned}       
\\&&  \hspace{-2em}
\mbox{where $\tau$ is assumed to be replaced by  $t_0+\delta t\,\sigma$}.
\nonumber
 \end{eqnarray}
They are polynomial in $\sigma,\mu(\tau),\lamBda(\tau)$.

We will sometimes call \Eq\eqref{int eq for mu} the first integral equation for \PIIItran{s}; accordingly,
\Eq\eqref{int eq for lamIII}
is the second integral equation.

The asserted interrelation of differential and integral equations in question is substantiated 
in appendices to the main text. %
In particular,
\Eqs\eqref{int eq for mu}, \eqref{int eq for lamIII} are derived from \Eqs\eqref{smoothness condition}, \eqref{ODE for mu with lamIII}, \eqref{ODE for lamIII}
in appendix \ref{app::A}.
For the sake of completeness,
a proof of fulfillment of \Eqs\eqref{smoothness condition}, \eqref{ODE for mu with lamIII}
in case of fulfillment of \Eq\eqref{int eq for mu}
is given in appendix \ref{app::B}.
Derivation of \Eq\eqref{ODE for lamIII} from  integral equations is more intricate since it requires preliminary establishing of appropriate regularity of the fraction
 \eqref{xi fraction} at the point $t=t_0$.
This can be carried out by means of introduction of the integral transformation sending the functions $ \mu(t),\lamBda(t) $ to
\begin{eqnarray}
                                            \label{upxi definition}              
&&
\begin{aligned}
\upxi(t)=&\;
-t_0^{-1}\big( (\chi_\infty +\sgn\,\chi_0 -1)/4 +2\mu(t) - 3 t_0 \,\lamBda(t) \big)
\\&\;
-t_0^{-1}\int^1_0 \!\!
d\sigma \, \sigma^3 \,\Omega_{\xi}[\mu,\lamBda](\sigma, t - t_0)
\end{aligned}
\end{eqnarray}
whose kernel is defined as follows
\begin{eqnarray}                        \label{Omega_xi definition}              
\hspace{-0em}
\Omega_{\xi}[\mu,\lamBda](\sigma,\delta t)
\!\!&=&\!\!
\sgn\,\big(
3(\chi_0-\sgn)
-8\chi_0\,\mu(\tau) - 3(\chi_0-2\sgn)\,t_0\,\lamBda(\tau)
)
\\ &&
   +3 \delta t \,\sigma\,
       t_0\,
       (2 \mu(\tau) - 1)
      \left(2  \sgn  \,\lamBda(\tau) + \big((\chi_0 -  \sgn )/(2t_0)  - \delta t \,\sigma\, \lamBda(\tau)\big)^2\right)
\nonumber\\ &&
   + 4 \delta t\, \sigma\,
       (\mu(\tau) - 1)\,\mu(\tau)\,
       \big(\sgn  -   \delta t\, \sigma\, (\chi_0 -  \sgn )/(2t_0) +  \sgn \, \delta t^2 \sigma^2 \lamBda(\tau)\big),
\nonumber
\\ && 
\hspace{-8.8em}
\mbox{where }\tau \mbox{ is assumed to be replaced everywhere with the expression } t_0+(t-t_0)\sigma.
\nonumber
\end{eqnarray}
The function $ \upxi(t) $ is determined by a non-local transformation 
and hence
is not, generally speaking, directly related to the function $\xi(t)$ defined point-wise by \Eq\eqref{xi fraction}
(it should be noted that, in spite of similarity, the notations $\upxi$ and $\xi$ still look different).
However, as it is shown in appendix \ref{app::C}, if
$\mu(t)$ and $\lamBda(t)$ are continuously differentiable
and the first integral equation \eqref{int eq for mu} is fulfilled then
 $\upxi(t)=\xi(t)$.

The above relation can be used for clarification of the degree of regularity
of the expression of right hand side of \Eq\eqref{ODE for lamIII}.
Specifically, in accordance with $\upxi(t)$ definition \eqref{upxi definition},
its smoothness degree is not smaller than the minimum of ones attributed to 
the functions
$\mu(t)$ and $\lamBda(t)$.
On the other hand,
it follows from definition
\eqref{xi fraction}  of $\xi(t)$ that in case of finiteness of the maximal degrees of smoothness of $\mu(t)$ 
and (or) $\lamBda(t)$ the maximal degree of its smoothness
at $t=t_0$ is one unit less.
In view of  coincidence of  $\xi(t)$  and $\upxi(t)$
taking place in case of fulfillment of \Eq\eqref{int eq for mu},
we must conclude,
avoiding contradiction, %
that
the smoothness of $\xi(t)$ 
is not limited from above. 

We come therefore to the following conclusions.
\begin{proposition}  \label{prop::070}
The equality of the expressions
defined by the formulas \eqref{xi fraction}  and \eqref{upxi definition}
which takes place 
for continuously differentiable solutions
$\mu,\lamBda$
to \Eq\eqref{int eq for mu}
implies
existence of derivatives of arbitrary order
for the fraction \eqref{xi fraction}.
\end{proposition}
\begin{corollary}  \label{coro::080}
If continuously differentiable functions  $\mu(t)$ and $\lamBda(t)$
obey the first integral equation \eqref{int eq for mu} then
the right hand side of
\Eq\eqref{ODE for lamIII}, irrespectively of its fulfillment,
admits derivatives of arbitrary order
everywhere including the point $t=t_0$.
\end{corollary}

The smoothness of the expression of the right hand side of the equation \eqref{ODE for lamIII}
on solutions to \Eq\eqref{int eq for mu}
allows further to prove its fulfillment on solutions to the system
of integral equations
\eqref{int eq for mu}
and
\eqref{int eq for lamIII}.
The details of the corresponding derivation can be found in appendix
\ref{app::D}.
Since the fulfillment of \Eq\eqref{ODE for mu with lamIII}
under compatible conditions, as well as inverse dependencies, has been established above,
we have the equivalence of the differential equations
with constrained initial date and the integral equations introduced above.
                    \begin{corollary}                    \label{coro::090}
On the class of continuously differentiable functions
the system of equations
 \eqref{smoothness condition}, \eqref{ODE for mu with lamIII},  \eqref{ODE for lamIII},
and
the system of equations
\eqref{int eq for mu},
\eqref{int eq for lamIII}
are equivalent.
                     \end{corollary}

It has to be noted that
there also exists yet another implementation of the second integral equation constituting in conjunction with
\Eq\eqref{int eq for mu} the system of equations equivalent to
\Eqs\eqref{int eq for mu}, \eqref{int eq for lamIII} and hence equivalent to
\Eqs\eqref{smoothness condition}, \eqref{ODE for mu with lamIII}, \eqref{ODE for lamIII}.
This modified ({\it alternative}) second integral equation looks as follows.
\begin{eqnarray}                  \label{another int eq for lamIII}               
&&\hspace{-2em}
\begin{aligned}[t]
\lamBda(t)=&\;  \lamBda(t_0)
-
{t-t_0\over t_0}
\left( \vphantom{ \int^1_0 }
(4t_0)^{-1}(\chi_\infty + \sgn\,\chi_0 - 1) + \lamBda(t)
 \right.\\[-1.5ex]&  \hspace{8.5em}\left.
-
(3t_0)^{-1} 
\int^1_0  \! d\sigma\:
\wideparen{\Omega}_{\lamBda}[\mu,\lamBda](\sigma,t-t_0)
\right),
\end{aligned}
\\  &&  \hspace{-2em}
\mbox{where }
                                    \label{parenOmega_lam definition}               
\begin{aligned}
\wideparen{\Omega}_{\lamBda}[\mu,\lamBda](\sigma,\delta t)
=&\;
2\Omega_{\mu}[\mu,\lamBda](\sigma,\delta t)
+
\sigma^3\,
\Omega_{\xi}[\mu,\lamBda](\sigma,\delta t).
\end{aligned}
\\[-1.5em]\hphantom{.}
\nonumber
\end{eqnarray}
It is worth mentioning that,
as opposed to \Eq\eqref{int eq for lamIII},
for $t=t_0$ \Eq\eqref{another int eq for lamIII} is trivially fulfilled
that simplifies analysis of approximate solutions and their convergence.

\begin{proposition}                   \label{prop::100}
The system of integral equations
\eqref{int eq for mu},
\eqref{int eq for lamIII}
is equivalent
of the system combining
the same first integral equation
\eqref{int eq for mu}
and the alternative second integral equation \eqref{another int eq for lamIII}.
\end{proposition}
\noindent
The asserted equivalence follows from computations set out in appendices
\ref{app::E}, \ref{app::F} to the main text and corollary \ref{coro::090}.

We will also name the system of integral equations \eqref{int eq for mu},
\eqref{another int eq for lamIII}
alternative
as against
the system of equations \eqref{int eq for mu},
\eqref{int eq for lamIII}.

\section{Analytic solutions to integral equations for \PainIII{ } transcendents and theorem \ref{teo::020}}\label{sec::040}

At many points of the above reasoning
 a proper regularity of the functions we consider was required. 
This precaution is justified 
since sometimes a survival of such a property at roots of a \PIIItran{} seems unobvious.
There are two cases which are here encountered.
Namely, basing
on differential equations
\eqref{Ham eq first}, \eqref{Ham eq second}, a relevant substantiation can be inferred from their regularity
everywhere except the center $t=0$.
On the other hand,  when starting with integral
equations, the regularity of their solutions should be demonstrated independently (although with the same restriction $t\not=0$).
Such a conclusion will be one of the outcomes of the discussion presented below.

\subsection{Iterative form of the alternative system of integral equation}

\Eqs{}\eqref{int eq for mu}, \eqref{another int eq for lamIII}
suit well for 
development of
an iterative algorithm enabling one
to construct  sequences of analytic functions (in fact, polynomials) $\mu_n(t),\lamBda_n(t),\,n=0,1,2,\dots$
which can be
regarded as their approximate
solutions
of successively increasing accuracy,
We consider below their minor modification allowing for a simple representation of such a procedure.

Generally speaking,
supposing
some starting element comprising the functions $\mu_0(t), \lamBda_0(t)  $ to be given,
subsequent elements are obtained by means of iterative application
of
the
pair of
subsequent
transformations
\begin{equation}                           \label{basic iteration}              
\begin{aligned}[t]
\mu_{n+1}(t) =&\;
 \raisebox{-1.2ex}[1ex]{\Big\lfloor}%
_{\hspace{-1.1ex}
 \mbox{\scriptsize
 \raisebox{0.9ex}[1ex]{
 \rlap{\:
 $
\mu
 \leftleftharpoons
 \mu_n,\;
 \lamBda
 \leftleftharpoons
 \lamBda_n
 $
 }
 }}
 }
{\sf RHS}[\mbox{Eq.\eqref{int eq for mu}}]
,
\;
\lamBda_{n+1}(t) =
 \raisebox{-1.2ex}[1ex]{\Big\lfloor}%
_{\hspace{-1.1ex}
 \mbox{\scriptsize
 \raisebox{0.9ex}[1ex]{
 \rlap{\:
 $
\mu
 \leftleftharpoons
 \mu_{n+1},\;
\lamBda
 \leftleftharpoons
 \lamBda_n
 $
 }
 }}
 }
{\sf RHS}[\mbox{Eq.\eqref{another int eq for lamIII}}].%
\end{aligned}
\end{equation}
It is obvious that if the functions
$ \mu_0(t), \lamBda_0(t)$ are smooth
then
the smoothness properties of all subsequent $\mu_n(t), \lamBda_n(t), \;    n\in\mathbb{N},$ are not worse, at least.

It is however more convenient to handle instead of the  sequence $\{\mu_n,\lamBda_n \} $
the associated sequence of
`the per step increments'
\begin{equation*} 
 \{ \deltamu_n, \delta\lamBda_n \}\!:\;
\deltamu_n(t)= \mu_n(t)-\mu_{n-1}(t)
,
\delta\lamBda_n(t)= \lamBda_n(t)-\lamBda_{n-1}(t)
.
\end{equation*} 
These differences  obey the equations looking like a chain of integral transformations
\begin{equation}
                                          \label{integral transforms}                          
\begin{aligned}[c]
{\deltamu_{n+1}}( t)
 =
 \hspace{3em}  & \hspace{-3em}
 -{{ \delta t}\over t_0}{\deltamu}_{n}( t)
 +{{ \delta t}\over t_0}
 \int_0^1 \!\! d\sigma\,
\delta\Omega_\mu[\mu_{n-1},\lamBda_{n-1},{\deltamu}_n,\delta\lamBda_n](\eta),
\\
{\delta\lamBda_{n+1}}( t)
 =
  \hspace{3em}  & \hspace{-3em}
 -{{ \delta t}\over t_0}{\delta\lamBda}_{n}(t)
-{{ \delta t}\over 3 t_0^2}
 \int_0^1 \!\!d\sigma
\,
\delta\wideparen{\Omega}_\lamBda[\mu_{n},\lamBda_{n-1},{\deltamu}_{n+1},\delta\lamBda_n](\sigma,\eta 
),
\end{aligned}
\end{equation}
where  $ \eta=\sigma\delta t, t=t_0+\delta t$, 
\\[0.5ex]
and the kernels
\begin{equation}                               \label{deltaOmegas decompositions}           
\hspace{-4.5em}
 \begin{aligned}[c]
\delta\Omega^{}_{\mu}[\mu,\lamBda,  \deltamu, \delta\lamBda ](\eta ) \hspace{2.4ex}  = &\;
\deltamu(\tau)\,
\delta\Omega^{(\mu)}_{\mu}[\hat\mu,\hat\lamBda](\eta)
+
\delta\lamBda(\tau)\,
\delta\Omega^{(\lamBda)}_{\mu}[\hat\mu,\deltamu](\eta )
,
\\
 \delta\wideparen{\Omega}_\lamBda [\mu,\lamBda,\deltamu,\delta\lamBda](\sigma,\eta)
= &\;
 2\delta\Omega_\mu[\mu,\lamBda,\deltamu,\delta\lamBda](\eta)
+\sigma^3\delta\Omega_\xi[\mu,\lamBda,\deltamu,\delta\lamBda](\eta)
\mbox{ with}
\\
\delta\Omega^{}_{\xi}[ \mu, \lamBda,\delta\lamBda, \deltamu ](\eta )  \hspace{2.4ex}   = &\;
\delta\lamBda(\tau)\,
\delta\Omega^{(\lamBda)}_{\xi}[\hat\mu,\hat\lamBda,\deltamu](\eta )
+
\deltamu(\tau)\,
\delta\Omega^{(\mu)}_{\xi}[\hat\mu,\hat\lamBda,  \delta\lamBda](\eta )
,
\end{aligned}
\end{equation}
\hspace{3ex}\\[-1.9em]
 \begin{equation*}
\hspace{-4.0em}
\mbox{where }
\tau=t_0+\eta,
\hat\mu(\tau)=\mu(\tau)+\deltamu(\tau)/2,
\hat\lamBda(\tau)=\lamBda(\tau)+\delta\lamBda(\tau)/2,
\mbox{ and}
\end{equation*}
\hspace{1ex}\hphantom{and}\\[-3.em]
%
%
%
\begin{eqnarray} %
\hspace{3.5em}
                                     \label{delta Omega^mu_mu definition}           
\llap{$\delta\Omega^{(\mu)}_{\mu}[\hat\mu,\hat\lamBda]$}(\eta ) \hspace{2.5ex}  &=& \!\! 
\;
\sgn\,\chi_0
-
2\eta\,
(2\hat\mu(\tau) - 1)
\left( \sgn - \eta\, \left(\chi_0 - \sgn\right)\!/(2t_0) + \eta^2\hat\lamBda(\tau) \right),
 \\
\hspace{2.3em}
                                     \label{delta Omega^lam_mu definition}          
\llap{$\delta\Omega^{(\lamBda)}_{\mu}[\hat\mu,\deltamu] $}(\eta ) \hspace{-1ex}  \hspace{2.5ex}  &=&\!\!
 - \half {\eta^3}
\left( \big(2\hat\mu(\tau) - 1\big)^2 - 1
+ \deltamu(\tau)^2
\right),
\\
\hspace{2em}
                                       \label{delta Omega^mu_xi definition}         
\llap{$\delta\Omega^{(\mu)}_{\xi}[\hat\mu,\hat\lamBda,\delta\lamBda]$}(\eta )      \!\!&=&\!\!
- 8 \sgn \,  \chi_0
 + 3  \eta \, ( \chi_0  - \sgn )^2/(2  t_0 )
\\
\hspace{2em}
&& %
 + 4  \eta\,
   \big ( 2 \hat\mu(\tau) - 1 + 3  t_0 \, \hat\lamBda(\tau) \big)\times
\nonumber \\[-1ex]  && \hphantom{+ 4  \eta\, \,}
    \big(\sgn  -  \eta \, ( \chi_0  - \sgn )/(2  t_0 ) +  \eta ^2 \hat\lamBda(\tau) \big)
\nonumber \\  &&
  -6  \eta ^3  t_0  \big( \hat\lamBda(\tau)^2 - \delta\lamBda(\tau)^2/4 \big),
\nonumber
\\
\hspace{2em}
                                      \label{delta Omega^lam_xi definition}         
\llap{$\delta\Omega^{(\lamBda)}_{\xi}[\hat\mu,\hat\lamBda,   \deltamu] $}(\eta ) \!\!&=&\!\!
   - 3  \sgn \,  t_0 \, (\chi_0 - 2  \sgn )
\\
\hspace{2em}
 &&
   +
     6  t_0 \, \eta \,
       (2 \hat\mu(\tau) - 1)\times
\nonumber\\[-1ex] && \hphantom{ + 6  t_0 \, \eta \,  }
       \big( \sgn  - \eta (\chi_0 -  \sgn )/(2  t_0 ) + \eta^2 \hat\lamBda(\tau)  \big)
\nonumber\\&&
    +
    4 \eta^3
       \big(
         \hat\mu(\tau) (\hat\mu(\tau) - 1) + \deltamu(\tau)^2/4
          \big)
\nonumber
\end{eqnarray}
do not depend on the enumerating index $n$.

If one assigns to the role of $\mu_0 , \lamBda_0  $ an exact solution $\mu, \lamBda $ to the system of equations
\eqref{int eq for mu},
\eqref{another int eq for lamIII} then the null sequence
$\deltamu_{n}(t)\equiv0\equiv\delta \lamBda_{n}(t) \,\forall\, n\ge0$
proves to serve the corresponding solution to \Eqs\eqref{integral transforms}.
Otherwise a sequence of nonzero functions arises and if they converge to zero fast enough and uniformly with respect to variation of
their argument then the sums
$\mu(t)=\mu(t_0)+\sum_0^\infty \deltamu_{n}(t),
\lamBda(t)=\lamBda(t_0)+\sum_0^\infty \delta\lamBda_{n}(t)
$ can be regarded as constituents of a solution to the equations which is searched for.

It is worth mentioning 
the following
specificity 
 of the iterative form of alternative system of integral equation for \PainIII{ }transcendents:
 its formulas \eqref{integral transforms} - \eqref{delta Omega^lam_xi definition} {\it do not contain\/} the parameter $\chi_\infty$.

\subsection{Initial case}

The development of an iterative scheme based of the transformations \eqref{integral transforms}
requires specification of appropriate initial state for the data processed.
Starting with arbitrary $\mu_{0}$, $\lamBda_{0}$, $\deltamu_{1}$, $ \delta\lamBda_{1} $,
some sequence of functions $\mu_{n-1}$, $ \lamBda_{n-1}$, $ \deltamu_{n}$, $ \delta\lamBda_{n}$
can be obtained
but it will not define approximate solutions to \Eqs\eqref{int eq for mu}, \eqref{another int eq for lamIII} which we intend to construct.
To obtain a proper approximation, it is necessary to ensure fulfillment
of this role already  by the first sequence element,
i.e.\ by the pairs
 $\{\mu_{0},\lamBda_{0}\}$ and $\{ \mu_{1}=\mu_{0}+\deltamu_{1}, \lamBda_{1}=\lamBda_{0}+\delta\lamBda_{1}\}$,
their accuracy orders being distinct by unit.

More concretely, the functions $\mu_{0},\lamBda_{0}$
and $\mu_{1}, \lamBda_{1}$ must be related by \Eqs\eqref{basic iteration}.
Such a preparation of the adapted initial data
for subsequent iterations
is carried out below. The corresponding explicit expressions can be considered as the necessary addendum to \Eqs\eqref{integral transforms}.
They have to be implemented at the initial stage of iterations considered below in section \ref{subs::50.30}.

However,
at the very beginning,
one may start from an arbitrary `approximation' even if it has apparently no relation to a solution we search for.
The most convenient choice of this kind is 
\begin{equation}
                                  \label{init}         
   \mu_{0}(t)\equiv0,\, \lamBda_{0}(t)\equiv0.
\end{equation}
Even then
a valid, though crudest, approximation is produced by means of application
to  $ \mu=\mu_{0}, \lamBda= \lamBda_{0}$
of the transformation
 \eqref{basic iteration} get with $n=0$.
Its result reads
\begin{equation} 
                            \label{setting mu_1}         
\hspace{-1.5ex}
\begin{aligned}
\hspace{-5.em}
\mu_1(t)=\delta \mu_1(t) =
 \hspace{0em}  & \hspace{-0em}
\,
\big(1-\sgn\,(\chi_0^2-1)/(2t_0)+ 3t_0\,\lamIII\big)/2
-{\mathring\delta{}t}\,({\chi_\infty+\sgn\chi_0 -1})/2,
\end{aligned}\hspace{-5em}
\end{equation}
\begin{equation}
                \label{setting lamIII_1}         
 \phantom{\;}\hspace{-5.3em}
 \begin{aligned}[t]
\lamBda_1(t)=\delta \lamBda_1(t) =
 \hspace{5em}  & \hspace{-5em}
\,
\lamIII
-{\mathring\delta{}t}\,\chi_\infty/(4t_0)
\\
&\hspace{-3.4em}
+{\mathring\delta{}t}^2
\sgn\,
\Big(1 -\left(\sgn\,\chi_0\left(\chi_0-\sgn\right)/t_0 - 3t_0\,\lamIII\right)\times
\\[-0.7ex]  &   \hspace{1.85em}
                               \left(\sgn\,(\chi_0^2-1)^{}/(2t_0) - 3t_0\,\lamIII\right)\!\!\Big)/10
\\[-1.0ex]
      &   \hspace{-3.4em}
+{\mathring\delta{}t}^2({\chi_\infty+\sgn\chi_0 -1})\times
\\
   &  \hspace{-1.2em}
\Big( \sgn \, \chi_0/(10t_0)
\\
   & \hspace{-1.2em}
- {\mathring\delta{}t}
\left((\chi_0-\sgn)^2/(4t_0) + \left(\chi_0^2 - 1\right)\!/(3t_0)    - 2\sgn\,t_0\,\lamIII\right)\!/6
\\
   & \hspace{-1.2em}
+ {\mathring\delta{}t}^2
(\chi_0-\sgn)
\left(\sgn\,(\chi_0^2-1)/(2t_0)  -  3t_0\,\lamIII\right)\!/28
   \Big)
\\
&   \hspace{-3.4em}
+{\mathring\delta{}t}^3
(\chi_0-\sgn)
\left( \left( \sgn\, (\chi_0^2-1)/(2t_0) - 3t_0\,\lamIII\right)^2 -1\right)\!/36
\\& \hspace{-3.4em}
-
{\mathring\delta{}t}^4
({\chi_\infty+\sgn\chi_0 -1})^2
\big(\sgn/7- {\mathring\delta{}t}\left(\chi_0 - \sgn\right)\!/20\big)/4,
\\  \mbox{{where }}            \hspace{9em}  & \hspace{-9em}  {\mathring\delta{}t}=(t-t_0)/t_0  \mbox{ and }\lamIII=\lamBda(t_0).
\end{aligned}\hspace{-5em}
\end{equation}
The functions
$ \mu_{1}(t),  \lamBda_{1}(t)$ are polynomials in deviation $t-t_0$ from $t_0$ of the degrees 1 and 5,
 respectively.
Their
``initial values'' evinced 
at $t_0$
agree with \Eq{}\eqref{smoothness condition}. It also holds
\begin{equation*}         
\lamBda_1(t_0)=\lamIII,
\end{equation*}
where 
the constant $\lamIII$ may be arbitrary.
This parameter, in conjunction with $\sgn$ (fixed in advance) which
selects the family of solutions and
is equal to either $1$ or $-1$, identifies the whole solution whose approximations are constructed.

The above quadruplet $\mu_{0},\lamBda_{0}, \delta\lamBda_{1},\deltamu_{1} $ can be used as the initial case for
an 
iterative procedure of generating approximate solutions to \Eqs\eqref{int eq for mu}, \eqref{another int eq for lamIII}
on the base of transformations \eqref{integral transforms}.
All they are polynomials in $t-t_0$ of growing degrees whose initial parts might be considered as partial sums of some  power series.
The following statement follows from definitions
\begin{proposition}                     \label{prop::110}
Let a sequence of pairs $\{\mu_n(t), \lamBda_n(t)\},\; n=0,1,\dots $ of continuous functions
uniformly converging to continuous functions $\mu(t), \lamBda(t)$
be given
such that
\Eqs\eqref{init}-\eqref{setting lamIII_1} holds true
and  \Eqs\eqref{integral transforms} are fulfilled for
$\deltamu_n(t)= \mu_n(t)-\mu_{n-1}(t)
,
\delta\lamBda_n(t)= \lamBda_n(t)-\lamBda_{n-1}(t),\; n=1,2,\dots$
Then the limiting functions $\mu(t), \lamBda(t)$ solve \Eqs\eqref{int eq for mu}, \eqref{another int eq for lamIII} as well.
\end{proposition}

It is worth noting that \Eqs\eqref{setting mu_1}, \eqref{setting lamIII_1} are just the only point where the parameter $\lambda_\infty$ appears.
The subsequent transformations of $\mu_1(t), \lamBda_1(t)$ prove, on their own, to be independent of it.

\subsection{Algorithm of iterative constructing of approximate solutions}\label{subs::50.30}

Equations \eqref{integral transforms} considered as a non-linear transformation
taking a pair of given functions
$\big\{\deltamu_{n}(t),\delta\lamBda_{n}(t)\big\}$ to a similar pair $\big\{\deltamu_{n+1}(t),\delta\lamBda_{n+1}(t)\big\}$,
are the main constituent of the iterative algorithm enabling one to construct a sequence of approximate solutions
to
\Eqs\eqref{int eq for mu}, \eqref{another int eq for lamIII} in the form of truncated power series.
More exactly, the corresponding iterative procedure is compiled from the following.

\medskip
\noindent
{\bf The input data} to be processed is defined as follows: 
             \begin{quote} 
let $n$ be a positive integer and
\\
let the four entire functions (e.g.\ polynomials) $\lamBda_{n-1}(t), \mu_{n-1}(t), \delta\lamBda_n(t),  \deltamu_n(t)$ be given.
\\[1ex]
\end{quote}
\vspace{-4ex}
{\bf The algorithm step} \label{thestep}                                                      
comprises the following   substeps:
\begin{quote}
\begin{enumerate}\rm
\item  
The functions $\lamBda_{n}(t), \mu_{n}(t)$ are constructed as follows:
\begin{equation}                           \label{substep1}                                     
\lamBda_{n}(t) = \lamBda_{n-1}(t) + \delta\lamBda_n(t)   , \;\mu_{n}(t) = \mu_{n-1}(t) + \deltamu_n(t) ,
\end{equation}
\item  
The function $ \deltamu_{n+1}(t) $ is constructed by means of the following integral transformation:
\begin{eqnarray}                          \label{substep2}                                     
&&\hspace{-3.5em}
\begin{aligned}
{\deltamu_{n+1}}( t)
 = \;
 \hspace{3em}  & \hspace{-3em}
\mathring\delta t
  \Big(
 -{\deltamu}_{n}( t)
 +
 \int_0^1 \!\! d\sigma\,
\big(%
\deltamu_n(\tau)\,
\delta\Omega^{(\mu)}_{\mu}[\mu_{n-1/2},\lamBda_{n-1/2}](\eta )
\\ &\hspace{5.25em}
+
\delta\lamBda_n(\tau)\,
\delta\Omega^{(\lamBda)}_{\mu}[\mu_{n-1/2},\deltamu_n](\eta )
\big)%
\Big),
\end{aligned}
\\
%
&&\hspace{-3.5em}
\begin{aligned}
\mbox{where } &
\mathring\delta t=\delta t/t_0=(t-t_0)/t_0,
\;
\tau=t_0+\eta, \; \eta=(t-t_0)\sigma,
\\ &\:
\mu_{n-1/2}(t)=\half\big( \mu_{n}(t)+\mu_{n-1}(t) \big) = \mu_{n-1}(t)+\deltamu_n(t )/2,
\\ &\;
\lamBda_{n-1/2}(t)=\half\big( \lamBda_{n}(t)+\lamBda_{n-1}(t) \big) = \lamBda_{n-1}(t)+\delta\lamBda_n(t )/2,
\end{aligned} \nonumber
\end{eqnarray}
and the kernel functions are defined in \Eqs\eqref{delta Omega^mu_mu definition}, \eqref{delta Omega^lam_mu definition}.
\item 
The function $ \mu_{n+1}(t) $ is constructed as follows:
\begin{equation}                        \label{substep3}                                     
\mu_{n+1}(t) = \mu_{n}(t) + \deltamu_{n+1}(t) .
\end{equation}
\item 
The function $ \delta\lamBda_{n+1}(t) $ is constructed by means of the following integral transformation:
\begin{eqnarray}                        \label{substep4}                                     
&&\hspace{-3.5em}
\begin{aligned}
{\delta\lamBda_{n+1}}( t)
 = \;
 \hspace{3em}  & \hspace{-3em}
\mathring\delta t
  \Big(
 -{\delta\lamBda}_{n}( t)
 -2(3t_0)^{-1} \!\!
\int_0^1 \!\! d\sigma\,
\big(%
\deltamu_{n+1}(\tau)\,
\delta\Omega^{(\mu)}_{\mu}[\mu_{n+1/2},\lamBda_{n-1/2}](\eta )
\\[-1ex]
   &\hspace{9.99em}
+
\delta\lamBda_n(\tau)\,
\delta\Omega^{(\lamBda)}_{\mu}[\mu_{n+1/2},\deltamu_{n+1}](\eta )
\big)
\\[-1ex]
   & \hspace{-1.75em}
-(3t_0)^{-1} \!\!
\int_0^1 \!\! d\sigma\,\sigma^3
\big(%
\deltamu_{n+1}(\tau)\,
\delta\Omega^{(\mu)}_{\xi}[\mu_{n+1/2},\lamBda_{n-1/2}, \delta\lamBda_n  ](\eta )
\\[-1.7ex]
  &\hspace{6.3em}
+
\delta\lamBda_n(\tau)\,
\delta\Omega^{(\lamBda)}_{\xi}[\mu_{n+1/2},\lamBda_{n-1/2}, \deltamu_{n+1} ](\eta )
\big)%
\Big),
\end{aligned}
\\
&&       \hspace{-0em}
\begin{aligned}[b]
\mbox{where } &
\mathring\delta t,
\;
\tau, \; \eta,\;\lamBda_{n-1/2}(t), \deltamu_{n+1}(t)
\mbox{ have been defined above,}
\\& \hspace{-2.8em}
\mbox{the kernel functions are defined in
\Eqs\eqref{delta Omega^mu_mu definition} - 
\eqref{delta Omega^lam_xi definition}, and}
\\ & \hspace{-2ex}
\mu_{n+1/2}(t)=\half\big( \mu_{n+1}(t)+\mu_{n}(t) \big) = \mu_{n}(t)+\deltamu_{n+1}(t )/2.
\end{aligned}\nonumber
\end{eqnarray}
\item 
The function $ \lamBda_{n+1}(t) $ is constructed as follows:
\begin{equation}                      \label{substep5}                            
\lamBda_{n+1}(t) = \lamBda_{n}(t) + \delta\lamBda_{n+1}(t) .
\end{equation}
\end{enumerate}
                      \end{quote} 
After the above substeps,
one gains 
the collection
 $ \mu_{n}(t)$, $\lamBda_{n}(t)$,  $  \deltamu_{n+1}(t)   $, $ \delta\lamBda_{n+1}(t)$
of the four functions
similar to (and interdigitated with) the initial quadruplet
and denoted by the the same 
kernel symbols
which are now 
endowed with enumerating indices incremented by unit.
These functions are also analytic and, in particular case of polynomial
input data, 
also polynomial.

\subsection{Sequence of approximate solutions to integral equations and its convergence}\label{subs::50.40}

Let us constrain the variation of $t$ to vicinity of the point $t=t_0\not=0$  assuming that 
\begin{equation}
                                            \label{what is alpha}                              
 \begin{aligned}
&|{\mathring\delta{}t}|  < \, \alpha \Leftrightarrow
|t-t_0| < |t_0|\,\alpha
\mbox{ for some positive constant }\alpha<1,
\end{aligned}
\end{equation}
and, in particular, ensuring fulfilment of the condition $t\not=0$.
We will represent this constraint by the inclusion relation $t \in U_\alpha$.
 Such a set $U_\alpha$ is open, simply connected, convex, with compact closure.

The  definitions \eqref{setting mu_1}, \eqref{setting lamIII_1}
of the functions $\mu_{1}(t),  \lamBda_{1}(t)$
enable one to determine
(if needed, by explicit though, in case of $\lamBda_{1}$, somewhat bulky expressions) the quantities
$M^{(\mu)},  M^{(\lamBda)}$ as the
  positive real numbers satisfying inequalities
\begin{equation}                          \label{Ms definitions}
\sup_{ t \in U_\alpha }  |\mu_1(t) | \le \half M^{(\mu)},\;
\sup_{ t \in U_\alpha }  |\lamBda_1(t) |\le \half M^{(\lamBda)}.
\end{equation} %
Let also $ B_\mu^{(\mu)}, B_\mu^{(\lamBda)}, B_\xi^{(\mu)}, B_\xi^{(\lamBda)}$
be the real positive numbers such that
\begin{equation} 
                                          \label{bounds for Omegas}                            
\sup_{} \big|\Omega_{\mbox{\tiny \ding{67} }}^{\mbox{\tiny(\ding{93}) }}\!\! [\cdots\!\,]( t)\big|
\le 
B_{\mbox{\tiny \ding{67} }}^{\mbox{\tiny(\ding{93}) }}\!\!
\end{equation} 
where  {\tiny\ding{67}}$\,\in\{\mu,\xi\}$ and {\scriptsize\ding{93}}$\,\in\{\mu,\lamBda\}$
denote  indices represented by one of the listed symbols, and
where $ [\cdots\!\,] $ stands for the list of functions serving arguments for the given $\Omega$-function
which
depends on the concretization of references
 {\tiny\ding{67}} and {\scriptsize\ding{93}}
and has to be get from the
corresponding definition of $\Omega$\,s among the ones
specified by \Eqs{}\eqref{delta Omega^mu_mu definition}-\eqref{delta Omega^lam_xi definition}.

The domain on which the suprema in inequalities \eqref{bounds for Omegas}
have
to be
determined 
is
as follows:
\begin{itemize}
\item all $ t \in U_\alpha $ are taken into account;
\item continuous functions $\hat\mu({t})$, $\hat\lamBda({t})$ are arbitrary ones obeying the inequalities
$
 |\hat\mu({t})|  \le  M^{(\mu)},
 |\hat\lamBda({t})|  \le  M^{(\lamBda)}
 $
for all $  {t} \in U_\alpha  $;
\item
continuous functions $\deltamu({t})$, $\delta\lamBda({t})$ are arbitrary ones obeying the inequalities
$ |\deltamu({t})|  \le  2M^{(\mu)}, |\delta\lamBda(t))| \le 2 M^{(\lamBda)} $
for all $  t \in U_\alpha  $.
\end{itemize}
%

\begin{remark}                                \label{rem::010}
\rm
The
above upper bounds
$  M^{ \mbox{\tiny(\ding{93})} } ,  B_{\mbox{\tiny \ding{67} }}^{\mbox{\tiny(\ding{93}) }}\!\! $
are clearly non-unique; in particular, they can unobstructedly be increased.
Their changes would also arise from variation of the parameter $\alpha$ affecting the domain $U_\alpha$.
Nevertheless, it is always possible to set them in such a way
that, considered as functions of $\alpha$ with other parameters fixed, they were  non-decreasing.
In the other words, it can be claimed that as $\alpha$ is decreasing (remaining positive), neither of them is increasing.
This means, in particular, that for any positive
$\alpha' <\alpha $ the `original' quantities
$  M^{ \mbox{\tiny(\ding{93})} } ,  B_{\mbox{\tiny \ding{67} }}^{\mbox{\tiny(\ding{93}) }}\!\! $
continue to
fulfill 
the above set of inequalities involving them 
in which
$\alpha$ is replaced with $\alpha'$ and the domain $U_\alpha$ is replaced with $U_{\alpha'}\ni t: |t-t_0|\le t_0 \,\alpha'$.
\end{remark}

Let now $N$ be a positive integer greater than one.
Let also
a
finite sequences of
quadruplets of
functions
$\lamBda_{n-1}(t), \mu_{n-1}(t), \delta\lamBda_n(t),  \deltamu_n(t)$
starting
with $  \mu_0, \lamBda_0, \deltamu_1, \delta\lamBda_1$
defined 
by \Eqs{}\eqref{init}, \eqref{setting mu_1}, \eqref{setting lamIII_1}
be given
for $ n=1,\dots,N  $.
Let they
conform pairwise to 
the 
algorithm specified above in section \ref{subs::50.30} for $ n=2,\dots,N  $.

It is easy to see that since
the transformations realizing a step of the algorithm take polynomials to polynomials
and since the initial quadruplet
$  \mu_0, \lamBda_0, \deltamu_1, \delta\lamBda_1$ consists of zeros and polynomials,
all the available functions 
endowed with greater index values
are polynomial in $t$ as well.
In particular, they are analytic everywhere. 

Finally, we require fulfillment of inequalities
\begin{equation}
                                                 \label{geometric progression as majorant}                             
|\deltamu_{n}(t)| <   \half M^{(\mu)} |{\beta\mathring\delta t}|^{n-1},
|\delta\lamBda_{n}(t)| < \half M^{(\lamBda)} |{\beta\mathring\delta t}|^{n-1},
\mbox{where }\mathring\delta t=(t-t_0)/t_0,
\end{equation}
with some real $\beta>1$  for all $ t\in U_\alpha $.
Let us note 
that
$\deltamu_1, \delta\lamBda_1$
conform to these constraints
for any $\beta\ge1$
in view of the very definitions \eqref{Ms definitions} of $ M^{(\mu)} $ and $ M^{(\lamBda)}$.
The set of finite sequences of polynomials satisfying the declared requirements is thus non-empty.
Notice also that, having some applicable $\beta$,
 it can be increased without violation of conditions \eqref{geometric progression as majorant}.

In view of 
interrelation \eqref{substep1} of the functions $\mu_{n}(t)$ and $\deltamu_{n}(t) $ 
and the vanishing of $\mu_0(t)$,
the
decompositions
$ \mu_{n}(t) = \sum_{k=1}^n \deltamu_k(t) $ take place.
Applying to it the first inequality 
\eqref{geometric progression as majorant},
one obtains 
\begin{equation}                                     \label{bound for mu_n}
|\mu_{n}(t) | <  \sum_{k'=0}^{n-1}\half M^{(\mu)} |{\beta\mathring\delta t}|^{k'}
=
\half M^{(\mu)}(1 -\beta|\mathring\delta t|)^{-1} (1 - |{\beta \mathring\delta t}|^{n}),\;n=1,\dots, N.           
\end{equation}

In accordance with definition of the adopted domain of variation of $t$ 
it holds $\beta | \mathring\delta t  | \le \beta\alpha $.
Getting, if necessary, a smaller $\alpha$, but keeping $M$-bounds unchanged,
we ensure fulfillment of the following
additional
constraint
\begin{equation}
                                                \label{ineq::470}
\beta\alpha \le 1/2.
\end{equation}
Then  for $t\in U_\alpha$ it holds
$\beta | \mathring\delta t  | < 1/2 $ and all the factors in the right hand side of \eqref{bound for mu_n}
are positive.
 In such a case this inequality and \eqref{ineq::470}  imply that
\begin{equation}                \label{another bound for mu_n}                    
|\mu_{n}(t) | <\half M^{(\mu)}(1 -\beta\alpha)^{-1}< M^{(\mu)}.
\end{equation}
As a consequence, we also have
\begin{equation}                \label{more bounds}                               
\begin{aligned}
|\mu_{n-1/2}(t)|= &\; \half|\mu_{n}(t) +\mu_{n-1}(t)|<M^{(\mu)},\;
\\
|\deltamu_{n}(t)|= &\; |\mu_{n}(t) - \mu_{n-1}(t)| \le |\mu_{n}(t)| + |\mu_{n-1}(t)|<2 M^{(\mu)}.
\end{aligned}
\end{equation}

Analogous inequalities
\begin{equation}                                  \label{bounds for lams}         
\begin{aligned}
|\lamBda_{n}(t) | < & \; M^{(\lamBda)},
|\lamBda_{n-1/2}(t)|  <M^{(\lamBda)},\;
|\delta\lamBda_{n}(t)|  <2 M^{(\lamBda)}
\end{aligned}
\end{equation}
are inferred
from \Eqs{}\eqref{substep1}, \eqref{substep5}, and the second inequality  \eqref{geometric progression as majorant}
in a similar way.

The uniform bounds from above imposed to the functions
$ \mu_{n-1/2},  \lamBda_{n-1/2}, \deltamu_{n}$ by   inequalities \eqref{more bounds}, \eqref{bounds for lams}
conform to the conditions of validity of inequalities \eqref{bounds for Omegas} 
in case of
$\Omega$-functions
get in the forms used in \Eq{}\eqref{substep2}.
Taking them into account
and applying  Cauchy-Bunyakovsky-Schwarz (in what follows,  CBS) inequality,
we obtain 
the bound from above
for  the right hand side of the noted equation
which then gives rise to the following inequality %
\begin{equation}
                                 \label{ineq::510}
 \begin{aligned}
|\deltamu_{n+1}(t)|
\le &\;
 |{\mathring\delta{}t}|
\Big(| \deltamu_{n}(t)|
+B_\mu^{(\mu)}
\lgrp\int_0^1  \!\!\! d\sigma
                   |\deltamu_{n}(\tau)|^2
\,\rgrp^{1/2}
\\ &  \hspace{5.6em}
+
B_\mu^{(\lamBda)}
\lgrp\int_0^1  \!\!\!  d\sigma
                   |\delta\lamBda_{n}(\tau)|^2
\,\rgrp^{1/2}
\Big).
\end{aligned}
\\[-0ex]
\end{equation}
Here, as above,
the replacing
$\tau\leftleftharpoons t_0+(t-t_0)\sigma$
has to be carried out prior to the integral evaluation. 

Upon the above replacement,  the
right hand sides
of  the inequalities \eqref{geometric progression as majorant}
 are modified by
the their multiplication  by $\sigma^{n-1}$. This means that the operator
$\int^1_0 d \sigma \times$,
when applied to the upper bounds for  $|\deltamu_{n}|^2$  and  $|\delta\lamBda_{n}|^2 $, 
integrates in fact the power function $\sigma^{2n-2}$ yielding the factor of $(2n-1)^{-1}$.
Its square root $(2n-1)^{-1/2}\le 1$ is then  incorporated into the resulting summands.
Thus, taking into account \eqref{geometric progression as majorant}, we obtain
the following consequences of the inequality \eqref{ineq::510}:
\begin{equation}                                             \label{ineq::520}
\hspace{-4em}
 \begin{aligned}[b]
|\deltamu_{n+1}(t)|
\le &\;
 { {\mathring\delta{}t} \over 2}
\Big(
\big(
1+ (2n-1)^{-{1\over2}} B_\mu^{(\mu)}
\big)
    M^{(\mu)}
+  (2n-1)^{-{1\over2}} B_\mu^{(\lamBda)}
   M^{(\lamBda)}
\Big) |\beta {\mathring\delta{}t}|^{n-1}
\\
<& \;
\half
\beta^{-1}
\Big(
(1+ B_\mu^{(\mu)})
    M^{(\mu)}
+ B_\mu^{(\lamBda)}
   M^{(\lamBda)}
\Big) |\beta {\mathring\delta{}t}|^{n}.
\end{aligned}\hspace{-2em}
\end{equation}

Now let us introduce 
the additional constraint
\begin{equation}                      \label{Q1 definition}   
 \beta\ge
Q_1:=
1+ B_\mu^{(\mu)}
+ B_\mu^{(\lamBda)}
   M^{(\lamBda)} \! / M^{(\mu)} 
\end{equation}
to admissible values of $\beta$ which, once adopted, holds true independently of the value of $n$.
Inequality \eqref{ineq::520} then implies that
\begin{equation}                      \label{reproducing of deltamu-bound}  
 \begin{aligned}
|\deltamu_{n+1}(t)|
< &\; \half M^{(\mu)}
|\beta {\mathring\delta{}t}|^{n}.
\end{aligned}
\end{equation}
Here $n=1,\dots,N$ and, thus, as compared to the original set  of such kind inequalities \eqref{geometric progression as majorant},
the next member indexed with $n=N+1$ has been added to it.
Besides, 
the next member of the sequence of functions
$\mu_n(t)$ defined by the formula \eqref{substep3}
obeys the inequality \eqref{bound for mu_n} and, in view of \eqref{ineq::470}, the inequality
\eqref{another bound for mu_n}  for $n=N+1$.

Inequalities \eqref{more bounds} with $n=N+1$ represent the additional upper bounds 
following from the above computations.

The relations derived above make it possible
to uniformly bound absolute values of $\Omega$-functions
from the integrals involved in \Eq{}\eqref{substep4}
by means of another application of inequalities \eqref{bounds for Omegas}.
Then applying to them CBS inequality
and replacing therein $\Omega$-functions by the (positive) constants $ B_{\mbox{\tiny \ding{67} }}^{\mbox{\tiny(\ding{93}) }}\!\! $,
one obtains
\begin{equation}                   \label{ineq::550}
 \begin{aligned}[b]
|\delta\lamBda_{n+1}(t)|
\le
 &\;
 |{\mathring\delta{}t}|
\bigg(
|\delta\lamBda_{n}(t)|
+2(3|t_0|)^{-1}
\Big(
 B_\mu^{(\mu)}
                  \lgrp\int_0^1  \!\!\!  d\sigma | \deltamu_{n+1}(\tau)  |^2\,\rgrp ^{1/2}
\\ & \hspace{10.26em}
+B_\mu^{(\lamBda)}
                  \lgrp\int_0^1  \!\!\!  d\sigma | \delta\lamBda_{n}(\tau)  |^2\,\rgrp^{1/2}
\Big)
\\
& \hspace{5.6em}
+(3|t_0|)^{-1}
\Big(
  B_\xi^{(\mu)}
                     \lgrp\int_0^1  \!\!\!  d\sigma \sigma^6 | \deltamu_{n+1}(t)  |^2\,\rgrp^{1/2}
\\ & \hspace{9.8em}
+ B_\xi^{(\lamBda)}
                     \lgrp\int_0^1  \!\!\!  d\sigma \sigma^6  | \delta\lamBda_{n}(t)  |^2\,\rgrp^{1/2}
\Big)
\bigg).
\end{aligned}
\end{equation}

In the same way as above,
inequalities \eqref{geometric progression as majorant}, \eqref{reproducing of deltamu-bound}
enable one
to reduce determination of the upper bounds for the residual integrals
in the right hand side  of
\eqref{ineq::550} 
to the integrating of the power functions $\sigma^{2n}$ and $\sigma^{2n-2}$  in the first  and second lines
and the power functions $\sigma^{2n+6}$ and $\sigma^{2n+4}$ in the third and fourth lines, respectively.
These lead to
 divisions of the coefficients associated with integrals by the square roots
$\sqrt{2n+1},\sqrt{2n-1}, \sqrt{2n+7},\sqrt{2n+5} $, respectively.
It follows therefore from
\eqref{ineq::550}  
that
\begin{equation}                        \label{ineq::560}
 \begin{aligned}[b]
|\delta\lamBda_{n+1}(t)|
\le
 &\;
 |{\mathring\delta{}t}|
\bigg(
\half M^{(\lamBda)}|\beta {\mathring\delta{}t}|^{n-1}
\\[-1ex]
 & \;\;\;\;\;\;
 +
 (3|t_0|)^{-1}
\lgrp 
{B_\mu^{(\mu)}M^{(\mu)}
\over
\sqrt{ 2n+1 }
}
|\beta {\mathring\delta{}t}|^{n}
+
{
  B_\mu^{(\lamBda)}  M^{(\lamBda)}\over\sqrt{2n-1}
}
|\beta {\mathring\delta{}t}|^{n-1}
\rgrp 
\end{aligned}\end{equation}\begin{equation}\nonumber\begin{aligned}
\\&\;\;\;\;\;\;
+(6|t_0|)^{-1}
\lgrp 
 {
   B_\xi^{(\mu)}M^{(\mu)}
\over
 \sqrt{2n+7}
}  |\beta {\mathring\delta{}t}|^{n}
+
 {
 B_\xi^{(\lamBda)} M^{(\lamBda)} \over\sqrt{2n+5}
}
|\beta {\mathring\delta{}t}|^{n-1}
\rgrp 
\bigg)
\\
 = &\;
(2\beta)^{-1}
\bigg(
\Big( 1
 +
{
  2 B_\mu^{(\lamBda)}  \over   3t_0     \sqrt{2n-1}
}
+
 {
 B_\xi^{(\lamBda)} \over    3t_0    \sqrt{2n+5}
 }
\Big)
 M^{(\lamBda)}
\\ &\hspace{4.5em}
+
{\beta |{\mathring\delta{}t}| \over   3t_0 }
\Big(
{
  2 B_\mu^{(\mu)}  \over       \sqrt{2n+1}
}
+
 {
 B_\xi^{(\mu)} \over        \sqrt{2n+7}
 }
\Big)
 M^{(\mu)}
\bigg)
|\beta {\mathring\delta{}t}|^{n}.
\end{aligned}
\end{equation}

Since $ |{\mathring\delta{}t}|< \alpha$ and in view of \eqref{ineq::470}
one has
$ \beta |{\mathring\delta{}t}| <1/2$.
For $n\ge 1$ the inequality \eqref{ineq::560} then implies that
\begin{equation}                \label{ineq::570}
 \begin{aligned}
|\delta\lamBda_{n+1}(t)| <
 &\;
(2\beta)^{-1}
\Big(
\big( 1
 +
(3t_0)^{-1}
(
{
  2 B_\mu^{(\lamBda)}
}
+
 {
 B_\xi^{(\lamBda)}  /       \sqrt{7}
 }
)
\big)
 M^{(\lamBda)}
\\ &\hspace{4.7em}
+
{(  6t_0)^{-1} }
\big(
{
  2 B_\mu^{(\mu)}  /       \sqrt{3}
}
+
 {
 B_\xi^{(\mu)} \!  /        3
 }
\big)
 M^{(\mu)}
\Big)|\beta {\mathring\delta{}t}|^{n}.
\end{aligned}
\end{equation}
We need to impose now to the auxiliary quantity $\beta$
yet another constraint which reads 
\begin{equation}
                                  \label{Q2 definition}  
 \begin{aligned}
\beta \ge
 Q_2 :=
& \; 1
 +
(3t_0)^{-1}
(
{
  2 B_\mu^{(\lamBda)}
}
+
 {
 B_\xi^{(\lamBda)}  /       \sqrt{7}
 }
)
 \\ & \hspace{0.8em}
+
{(  6t_0)^{-1} }
\big(
{
  2 B_\mu^{(\mu)}  /       \sqrt{3}
}
+
 {
 B_\xi^{(\mu)} \! /        3
 }
\big)
 M^{(\mu)} \!  / M^{(\lamBda)}.
\end{aligned}
\end{equation}
Once adopted, this relation holds true independently of the value of $n$.
With such $\beta$
the inequalities
\eqref{ineq::570}
give rise to the weaker ones which read
\begin{equation}                                \label{reproducing of deltalamIII-bound} 
 \begin{aligned}
|\delta\lamBda_{n+1}(t)|
\le
 &\;
\half M^{(\lamBda)}|\beta {\mathring\delta{}t}|^{n}.
\end{aligned}
\end{equation}
In the case of $n=N$ it
represents
the next member of the set of similar inequalities
\eqref{geometric progression as majorant} assumed to be fulfilled in cases of smaller $n$.

\Eqs{}\eqref{init},\eqref{substep5} imply the existence
of representation of
$\lamBda_{n+1}$
in the form of  a
 truncated series $ \lamBda_{n+1}(t)= \sum_{k=1}^{n+1} \delta\lamBda_{k}(t)  $.
Together with inequalities 
\eqref{geometric progression as majorant},  \eqref{ineq::470}, \eqref{reproducing of deltalamIII-bound}
this decomposition
leads to the inequalities
\begin{equation}                                   \label{reproducing of lamIII-bound} 
 \begin{aligned}
|\lamBda_{n+1}(t) | \le & \; 
\half M^{(\lamBda)}(1 -\beta|\mathring\delta t|)^{-1} (1 - |{\beta \mathring\delta t}|^{n+1})
<\half M^{(\lamBda)}(1 -\beta|\mathring\delta t|)^{-1}
\\ < & \;
\half M^{(\lamBda)}(1 -\beta\alpha)^{-1}
< M^{(\lamBda)}.
\end{aligned}
\end{equation}
In case of
$n=N$
this bound
is equivalent to
the next member of the
first sequence of inequalities \eqref{bounds for lams}. The elongation 
by one element of the second sequence from the triplet 
then appears
as well.

The above estimates would establish the concordance of the algorithm
of
constructing of
functions
$\mu_{n}(t), \lamBda_{n}(t), \deltamu_{n+1}(t), \delta\lamBda_{n+1}(t)  $
from the functions
$\mu_{n-1}(t),$ $ \lamBda_{n-1}(t),  \deltamu_n(t), \delta\lamBda_n(t) $
formulated in section \ref{subs::50.30}
to the properties which they are assumed to be endowed with. 
However,  one point concerning their capability
still needs clarification.
It concerns  the consistency
of the restrictions
imposed
to the auxiliary parameters
$\alpha$ and $\beta$.

As to $\beta$, it has to obey inequalities \eqref{Q1 definition}, \eqref{Q2 definition}.
Hence the condition
\begin{equation}                       \label{min for beta}  
\beta \ge \max(Q_1,Q_2) > 1
\end{equation}
is to be met. However, the inequality \eqref{ineq::470} may prove to be inconsistent with it.

Fortunately,
such a fairly
predictable quandary can be prevented by means of a small modification of the foregoing reasoning.
Its outline
constitutes a part of the proof of the following
                          \begin{theorem} \label{teo::120}
Let the functions $\mu_1(t), \lamBda_1(t)$ 
be defined by the formulas \eqref{setting mu_1},\eqref{setting lamIII_1},
where $
{\mathring\delta{}t}$ denotes `the normalized deviation' $(t-t_0)/t_0$ of the argument $t$ from the given point $t_0\not=0$;
\\
--- let the positive numbers $M^{(\mu)}, M^{(\lamBda)} $ be defined
as follows:
\begin{equation}
                                                \label{eq::620}
M^{(\mu)} = \max(1,2\sup_{t \in U_1 } |\mu_1(t)|),
\;
M^{(\lamBda)}=\max(1,2
\sup_{t \in U_1 }
|\lamBda_1(t)|),
\end{equation}
where $U_1$ is the open disc centered at $t_0$ and touching  zero, $U_1\ni t  \Leftrightarrow |\mathring\delta{}t| <1 \Leftrightarrow |t-t_0|<|t_0|$;
\\
--- let the positive numbers
$B_\mu^{(\mu)}, B_\mu^{(\lamBda)},  B_\xi^{(\mu)}, B_\xi^{(\lamBda)} $
be defined
by the formulas
\begin{equation}
                                               \label{eq::630}
B_{\mbox{\tiny \ding{67} }}^{(\mbox{\tiny\ding{93}})}
=\max(1,\sup_{t \in U_1 } |\Omega_{\mbox{\tiny \ding{67} }}^{(\mbox{\tiny\ding{93}})}\! [\cdots\,\!]( t)|),
\end{equation}
where the symbols
{\scriptsize\ding{67}}$\,\in\{\mu,\xi\}$ and {\scriptsize\ding{93}}$\,\in\{\mu,\lamBda  \}$
stand for the symbolic `indices' 
which label the four functions
$ \Omega_{\mbox{\tiny \ding{67} }}^{(\mbox{\tiny\ding{93}})}  $
defined by the formulas \eqref{delta Omega^mu_mu definition}-\eqref{delta Omega^lam_xi definition}.
In them,
the lists in square brackets
imaged in \eqref{eq::630} as   $[\cdots\,\!]$
contain
`the functional arguments'
which refer to
continuous functions obeying inequalities
$$
 |\hat\mu(t)|  \le  M^{(\mu)}  ,
  |\hat\lamBda(t)|  \le  M^{(\lamBda)} ,
 |\deltamu(t)| \le 2 M^{(\mu)}  ,
  |\delta\lamBda(t)\le  2M^{(\lamBda)}
$$
for all $t \in U_1 $;
\\
--- let the positive numbers $Q_1,Q_2$ are defined by the formulas
\eqref{Q1 definition} and  \eqref{Q2 definition}, respectively;
\\
--- let the positive numbers $\beta,\tilde\alpha $ be defined as follows:
\begin{equation}                    \label{beta def}                              
\beta = \max(Q_1,Q_2),
\;
\tilde\alpha = (2\beta)^{-1};
\end{equation}
--- let the sequences
 $
\mu_n(t),
 \lamBda_n(t),
\deltamu_n(t),
\delta\lamBda_n(t),
\; n=1,2,\dots$
of polynomials in $t$
be constructed
using the iterative transformations 
\eqref{substep1}---\eqref{substep5}
and
starting from  $\mu_1(t), \lamBda_1(t)$ specified by the formulas \eqref{setting mu_1}, \eqref{setting lamIII_1}
and $ \delta \mu_1(t)=\mu_1(t), \delta\lamBda_1(t)=\lamBda_1(t)$.

Then the polynomial sequences $\{\mu_n(t)\}  ,\{\lamBda_n(t)\}$ uniformly converge on  the
domain  $ U_{\tilde\alpha}\ni t: |t-t_0| < |t_0|\tilde\alpha$
to analytic functions $\mu(t), \lamBda(t)$, respectively. The latter verify the integral equations
\eqref{int eq for mu}, \eqref{another int eq for lamIII}
and obey the inequalities
\begin{equation}                                  \label{bound for mu}                             
 \begin{aligned}
|\mu(t) | \le & \;
\half M^{(\mu)}(1 -\beta|t_0|^{-1} |t-t_0|)^{-1},
\end{aligned}
\end{equation}
\begin{equation}                                   \label{bound for lam}                           
 \begin{aligned}
|\lamBda(t) | \le & \;
\half M^{(\lamBda)}(1 -\beta|t_0|^{-1} |t-t_0|)^{-1}.
\end{aligned}
\end{equation}
                       \end{theorem}
                       \begin{proof}
The definitions \eqref{eq::620}, \eqref{eq::630} agree with   conditions \eqref{Ms definitions}, \eqref{bounds for Omegas} considered in case $\alpha=1$.
Intending to apply the mathematical induction, we assume the existence of a finite sequence of quadruplets of functions
$\mu_{n-1}, \lamBda_{n-1}, \deltamu_{n}, \delta\lamBda_{n},\; n=1\dots,N \ge 1$ (in which $\mu_0(t)\equiv0\equiv\lamBda_0(t)$) obeying inequalities
\eqref{geometric progression as majorant}. Since in case $n=1$ the latter are fulfilled by virtue of definitions \eqref{eq::620},
the set of such sequences is not void.

Inequalities
\eqref{geometric progression as majorant}
lead to validity of the upper bounds \eqref{bound for mu_n} taking place for $n=1,\dots,N$.
We will consider them on the smaller domain $U_{\tilde\alpha}\subset U_1$. For all $t\in U_{\tilde\alpha}$
it holds $|\beta \mathring\delta t|<1/2$ and inequalities \eqref{bound for mu_n} then imply existence of the upper bound \eqref{another bound for mu_n}
and subsidiary relations \eqref{more bounds}.

Now all the conditions which had led to inequality \eqref{ineq::520} are fulfilled.
Since inequality \eqref{Q1 definition} takes place in accordance with theorem conditions, the existence of the upper bound
\eqref{reproducing of deltamu-bound} follows.

Inequalities
\eqref{bounds for lams}, \eqref{reproducing of deltalamIII-bound}, \eqref{reproducing of lamIII-bound}
are then obtainable by means of repetition of their original derivation.

Index $n=N$ in the above inequalities corresponds to the functions (in fact, polynomials)
$\mu_N,\lamBda_N,\deltamu_{N+1},\delta\lamBda_{N+1}$.
We have proven therefore the validity  of the supposed upper bounds
for the next element of the known finite non-empty sequence of function quadruplets pairwise bound in accordance with the algorithm specified
in section \ref{subs::50.30} on the domain $U_{\tilde\alpha}$.
Then the principle of mathematical induction establishes existence of the infinite sequence containing and extending all the above finite ones. 

We have the expansions
$\mu_n(t)=\sum_{k=1}^{n-1} \deltamu_{k}(t)$,
$\lamBda_n(t)=\sum_{k=1}^{n-1} \delta\lamBda_{k}(t)$
which in conjunction with inequalities \eqref{geometric progression as majorant}
prove the convergence of these series to analytic functions
$\mu(t)=\lim_{n\to\infty} \mu_{n}(t)$,
$\lamBda(t)=\lim_{n\to\infty} \lamBda_{n}(t)$,
provided $ |\mathring\delta t|<\beta^{-1}$  or, equivalently, $|t-t_0|< |t_0|/\beta$.
The choice of the domain $U_{\tilde\alpha}$ of $t$ variation verifies the latter condition.

The limits for the function sequences $\deltamu_n(t), \delta\lamBda_n(t)$ as $n\to\infty$ also exist. They coincide with the identically
zero functions to which the convergence on $U_{\tilde\alpha}$ is uniform. Since
$\deltamu_n(t)= \mu_n(t)-\mu_{n-1}(t), \delta\lamBda_n(t)= \lamBda_n(t)-\lamBda_{n-1}(t)$
these functions with index $n$ differ from functions with index $n-1$ by the quantities uniformly tending to zero.
We can pass to the limit as $n\to\infty$ in the equations \eqref{basic iteration}
which are equivalent to the relations \eqref{substep2}, \eqref{substep4} and hence are fulfilled.
The result arises after a mere discarding of indices marking functions with kernel names  $\mu$ and $\lamBda$.
It coincides with \Eqs{}\eqref{int eq for mu}, \eqref{another int eq for lamIII}
which are therefore fulfilled.

Let us note in conclusion of the proof that
the bounding
\eqref{bound for mu} follows  from inequalities $ |\mu_{n}(t) | <\half M^{(\mu)}(1 -\beta\alpha)^{-1} $
which are a part of inequality sequence \eqref{another bound for mu_n}.
Similarly, inequality
\eqref{bound for lam} follows from the sequence of inequalities $ |\lamBda_{n+1}(t) | \le  \half M^{(\lamBda)}(1 -\beta|\mathring\delta t|)^{-1} $
embedded in the chain \eqref{reproducing of lamIII-bound}.
\end{proof}
Having analytic solution to integral equations \eqref{int eq for mu}, \eqref{another int eq for lamIII},
one has in fact an analytic solution to differential equations \eqref{ODE for mu with lamIII}, \eqref{ODE for lamIII}.
Then the formula \eqref{eq::120} defines the analytic function $\lambda(t)$ which verifies
\Eq\eqref{equation PIII'} and vanishes at $t_0$.
The characteristics of the set of such solutions
asserted in theorem \ref{teo::020} follow from the form of the iteration
starting data \eqref{setting mu_1}, \eqref{setting lamIII_1}
taking into account \Eq\eqref{this is sign switch}
and the same decomposition  \eqref{eq::120}.
These relations prove theorem \ref{teo::020} which is thus a consequence of theorem \ref{teo::120}.\hfill$\square$

\section{Bounding minimal distance between a non-zero root and other roots and poles}\label{sec::050}

In accordance with one of definitions, a \PIIItran{} possesses the Painlev\'e property
allowing only to poles to be situated at the points   of its non-analyticity.
The only exception is the possible irregularity at zero where \Eq\eqref{equation PIII'} is singular itself.
At the same time the theorem \ref{teo::120} states that $\lambda(t)$ is analytic inside the disk of the radius $\tilde\alpha\,t_0$  centered at $t_0$.
Hence all the poles, if there are any, are situated outside the latter.
This relation yields a quantitative bound from below for the distance between a non-zero
root and a pole of $\lambda(t)$. Namely, we have the following theorem consequence.
               \begin{corollary}\label{coro::130}
Let $t_0\not=0$ be a root of a \PainIII{ }transcendent $\lambda(t)$.  Then the distance between $t_0$ and any pole of $\lambda$
is not less than $\tilde\alpha\,|t_0|$, where $\tilde\alpha$ is defined by \Eqs\eqref{beta def}.
               \end{corollary}

Associated results can also be used for deriving an estimate from below for the distance between the roots themselves.
To that end, let us note that \Eq\eqref{eq::120}
rewritten as follow
\begin{equation*}\nonumber
\lambda(t) = (t-t_0)\left(
\sgn + \half t_0^{-1}(\sgn - \chi_0)(t-t_0) + (t-t_0)^2 \lamBda(t)
\right)
\end{equation*}
implies in view of \Eq\eqref{this is sign switch} the fulfillment of the inequality
\begin{equation*}\nonumber
|\lambda(t)|
\ge
|t-t_0|\,
\left(
1 - |\half t_0^{-1}(\sgn - \chi_0)|\,|t-t_0| - |t-t_0|^2 |\lamBda(t)|
\right)
 \end{equation*}
 Using the inequality \eqref{bound for lam}
%
one further obtains that
\begin{equation*}
|\lambda(t)|
>
|t-t_0|  
\left(
1 - |\half t_0^{-1}(\sgn - \chi_0)|\,|t-t_0|
- {
\half M^{(\lamBda)}
|t-t_0|^2
\over
1 -\beta|t_0|^{-1}|t-t_0|
}
\right)
\end{equation*}
leading finally to the inequality
\begin{eqnarray}                     \label{ineq::670}
|\lambda(t)|
&> &
{|\delta t|
\over
1 -(2\tilde\alpha\,|t_0|)^{-1} |\delta t|
}\: Q(|\delta t|), \mbox{ where } \delta t=t-t_0 
\\
\llap{\mbox{and }}
Q(x)  \!\! &=& \!\!
\big(1 - |\half t_0^{-1}(\sgn - \chi_0)|\,x\big)\big(1 -(2\tilde\alpha\,|t_0|)^{-1} x\big)
- {
\half M^{(\lamBda)}
x^2
}.
\nonumber
\end{eqnarray}
Concordantly with the role of $t_0$,
the right hand side of \eqref{ineq::670} 
vanishes as $\delta t=0$. 
Its other vanishing may occur either at a root of the quadratic function $Q$ or sufficiently far of $t_0$ where
the substantiation of \eqref{ineq::670}
ceases to be proper
and its
fulfillment
is not guarantied.
Accordingly, the following statement holds true.
                  \begin{theorem}\label{teo::140}
Let the quadratic equation $Q(x)=0$ has real roots and let the minimal positive one among them, $x_{min}$,
obeys the condition $x_{min}\le\tilde\alpha$. 
Then all the nonzero roots of \PainIII{} transcendent $\lambda(t)$ distinct of $t_0$
are away of $t_0$ at distances not smaller than $x_{min}|t_0|$.
Otherwise, such distances cannot be smaller than $\tilde\alpha|t_0|$.
                  \end{theorem}
The quantities appearing the above estimates can be computed by means of explicit formulas given in the preceding section.
\begin{remark}          \label{rem::020}
\rm
For each pairs of roots of a \PIIItran, interchanging them, two bounds from below for their mutual distance arises.
Obviously, the maximal one among them   has to be get as the result.
\end{remark}

\section{Approximate representation of \PainIII{} transcendent near its root}\label{sec::060}

\subsection{Derivation of explicit approximate solution} 

The systems of two coupled integral equations
\eqref{int eq for mu},
\eqref{int eq for lamIII}
which an appropriately transformed \PIIItran{} obeys
can be used as the base for   iterative computational scheme
enabling one to construct its explicit approximate representation in terms of truncated power series
applicable in vicinities of zeros.
As compared to the alternative second integral equation
\Eq\eqref{another int eq for lamIII}, using of \Eq\eqref{int eq for lamIII}
is here preferable since the kernel of the integral transformation in the former is more complicated. 

The sought approximations to a \PIIItran{}
arise as elements of
the sequences of polynomials
$\mu_m(t), \lamBda_n(t)\: (m,n\in\mathbb{Z}),$ 
determined from 
an appropriate subset selected from families of the equations
\begin{eqnarray}
                             \label{iterative rule for mu}         
&& \hspace{-3em}
\begin{aligned}
\mu_{m_1+1}(t_0+\delta t)=&\;
\half
\big(
    1 - \sgn\,(\chi_0^2 - 1)/(2t_0) +3\,t_0\,  \lamIII
\big)
\\
&\;
+
  \delta  t\,t_0^{-1}
  \big(
       - \half(\chi_\infty+\sgn\,\chi_0-1) - \mu_{m_1}(t_0+\delta t)
\\&\; \hphantom{+ \delta  t\,t_0^{-1}  \big(   } \;
  +
   \int^1_0 \!\!
  d \sigma\:
  \Omega_{\mu}[\mu_{m_1},\lamBda_{m_2}](\sigma,\delta  t) \big),
\end{aligned}
\\
                                                   \label{iterative rule for lamIII}           
&& \hspace{-3em}
\begin{aligned}
\lamBda_{n_1+1}(t_0+\delta t)
=&
-
\delta  t\,
t_0^{-1}
\lamBda_{n_1}(t_0+\delta t)
+
t_0^{-1}\!\!
\int_0^1 \!\! d \sigma\,
\sigma^2\,
\Omega_{\lamBda}[\mu_{n_2},\lamBda_{n_1}](\sigma,\delta  t),
\end{aligned}
\end{eqnarray}
where either $\sgn=1, $ or $\sgn=-1, $ but not both,
 $\lamIII=$\;const is arbitrary fixed,
and $m_1,m_2,n_1,n_2\in{\mathbb Z}$.
These are here considered as transformations converting the pair of
functions
 $ \mu_{m_1},\lamBda_{m_2} $ to the function $ \mu_{m_1+1}$
and the pair
$\mu_{n_2},\lamBda_{n_1} $ to $ \lamBda_{n_1+1}$.

It is easy to see that if the above sequences are right-side infinite and if  $\mu_m(t)$ and $\lamBda_n(t)$ converge
uniformly on some domain with compact closure containing $t_0$ then their limits  $\mu(t)$, $\lamBda(t)$
verify \Eqs\eqref{int eq for mu}, \eqref{int eq for lamIII}
and thus define a \PIIItran{ }vanishing at $t_0$.
At the same time separate sequence elements
can play role of approximate solutions
to the same equations.
We will see that they are polynomial in  $t$ or, equivalently but more appropriately, in $\delta t =t-t_0$,
and that their parts of lower degrees a kept unchanged in the subsequent elements of these sequences.
Moreover, they are not affected if with each step to the right along the sequence the terms of degrees greater than some threshold are dropped out
although  then, strictly speaking, the fulfillment of \Eqs\eqref{iterative rule for mu}, \eqref{iterative rule for lamIII}
becomes violated.

To proceed, it is necessary to choose the form of the 
initial case allowing one to efficiently carry out subsequent computations leading to a meaningful result.
The general observation seeming here relevant is as follows:
we should agree that, at the beginning, we know nothing  about the solution we plan to obtain.
Accordingly, as above, we may simply assign
\begin{equation}
                                                         \label{eq::700}
\mu_{0}(t)\equiv \lamBda_{0}(t)\equiv 0.
\end{equation}
It is natural to regard this `approximation' as
containing {\it no\/} genuine information about the
limiting
functions for the sequence conformal to the relations \eqref{iterative rule for mu}, \eqref{iterative rule for lamIII}.

Nevertheless,  one may 
 substitute zeros \eqref{eq::700} 
into the right hand sides of the noted formulas get with $m_1=m_2=n_1=n_2=0$.
The result is
\begin{equation*}\nonumber
\mu_{1}(t)=\half\big(1 -(\chi_0^2-1)/(2t_0) +3 t_0\,\lamIII\big),\;
 \lamBda_{1}(t)=\lamIII.
\end{equation*}

It has to be noted that it does not agree with
\Eqs\eqref{setting mu_1},\eqref{setting lamIII_1}.
Nevertheless,
the above formulas prove to already gain some relevant data.
Although these functions  reduce to constants, the latter coincide with
the values of the corresponding limiting functions (exact solutions) at the point $t=t_0$.
Hence they can be considered as the implementation of the ``zero order approximation''
(being labeled by the indices 1, though)
which will be ``embedded'' unchanged in all the subsequent higher order ones.

But at the moment, we still know nothing about how $\mu$ and $\lamBda$ vary under the impact of
variation of $t$.
In response to this, let us recast the above equalities as follows:
\begin{equation}                        \label{eq::710}
\mu_{1}(t_0+\delta t)=\half\big(1 -(\chi_0^2-1)/(2t_0) +3 t_0\,\lamIII\big) + \Or(\delta t),\;
 \lamBda_{1}(t_0+\delta t)=\lamIII + \Or(\delta t).
\end{equation}
Here  $\Or$-terms
can be regarded as existing but currently unknown contributions vanishing at $t_0$.
Generally speaking, contributions of the equal or higher orders generated by any algebraic operation (summation, multiplication) involving a quantity
declared unknown have to be declared unknown as well.

It is now pertinent to mention that the  kernels $\Omega_\mu, \Omega_\lamBda$  in integrals in formulas
\eqref{iterative rule for mu}, \eqref{iterative rule for lamIII}
 are polynomials  in  $\mu$ and $\lamBda$ --- as well as in all other explicit and suppressed arguments except of the constant $t_0$
which appears somewhere as a multiplier in denominators of fractions. See definitions \eqref{Omega_mu definition}, \eqref{Omega_lam definition}.
This means that if $\mu_*(t)$ and $\lamBda_{\star}(t)$ (where values of indices $*,\star\in\mathbb{N}$) are polynomial in $t$ or, equivalently, in $\delta t$,
then the integrands in formulas \eqref{iterative rule for mu}, \eqref{iterative rule for lamIII} are also polynomial in both $\delta t$ and $\sigma$.
Accordingly, the integrating over $\sigma$
proves to be
elementary
and, besides, the polynomial dependence on $\delta t$ is retained. 
The following conclusion recapitulates the above notes.
                           \begin{proposition}\label{prop::150}
If the functions $\mu_{m_1},\lamBda_{m_2},\mu_{n_2},\lamBda_{n_1}\;  (m_1,m_2,n_1,n_2\in{\mathbb N})$ are polynomial in $t$
(or, equivalently, in $\delta t$)
then the functions
$\mu_{m_1+1},\lamBda_{n_1+1}$
obtained from them by means of the transformations
\eqref{iterative rule for mu},
\eqref{iterative rule for lamIII}
are also polynomial in $t$ (in $\delta t$).
                            \end{proposition}

Being devoid of their ``unknown constituents'', the functions \eqref{eq::710}
represent particular instances of polynomials.
Accordingly, the
iterative transformations \eqref{iterative rule for mu}, \eqref{iterative rule for lamIII},
employing them as the data enabling to start,
will also produce in all orders polynomials in $t$ or, equivalently, in $\delta t$ interpretable as truncated power series.
Elimination of terms ``of superfluous accuracy'' proving to be mixed with ``unknown contributions'' would retain such an interpretation.

The following useful observation concerns the important property of the kernel function $\Omega_\mu$
implied by its definition \eqref{Omega_mu definition}.
It shows explicitly the distinction of `the degrees of influence' of the functions $\mu$ and $\lamBda$ to the kernel value.
Namely,
the function $\mu$ is situated there, above all else, as the common factor affecting that value without any suppression.
On the contrary,
the magnitude of $\lamBda$
is suppressed by the factor of $\delta t^3$ of the third order of smallness.

Hence,
when computing the functions $\mu_{m+1}$ with the required accuracy $m$ (i.e.\ up to an `unknown' contribution of the degrees $m+1$ and higher)
by means of the transformation \eqref{iterative rule for mu}, one must use on the right
$\mu$-function given with accuracy of the degree $m-1$, i.e.~$\mu_{m}$
(an additional unit in the degree is ensured by the factor of $\delta t$ situated in the second line
of~\eqref{iterative rule for mu}).
At the same time, it is enough for $\lamBda$-functions involved in the same calculation
 to be three orders less accurate, i.e.~with the degree equal to $m-4$ --- but not less than that.
The terms with degrees of $\delta t$ higher than $m-4$ are here irrelevant since their inclusion would affect only `uknown' contribution
remaining as such
and hence
would have no impact on
the meaningful part of the final result.

The above remarks apply to the following initial steps of computations of approximate solutions to
the corresponding equations \eqref{iterative rule for mu}, \eqref{iterative rule for lamIII}:
\begin{equation}                             \label{eq::720}
\{\mu_1,\lamBda_1\} \mbox{\,{\scriptsize$\models$}\:}  \mu_2,
\;\{\mu_2,\lamBda_1\} \mbox{\,{\scriptsize$\models$}\:}  \mu_3,
\;\{\mu_3,\lamBda_1\}  \mbox{\,{\scriptsize$\models$}\:} \mu_4,
 \;\{\mu_4,\lamBda_1\}  \mbox{\,{\scriptsize$\models$}\:} \mu_5.
\end{equation}

Here the symbol `{\scriptsize$\models$}' denotes the map applying the formula
\eqref{iterative rule for mu} with corresponding values $m_1=1,2,3,4$ and $m_2=1$
to functions between braces
which does not output, however, the result `as is' but interprets the `input' functional arguments
as combinations incorporating certain `known' part and some `unknown' one, the latter being of a higher degree in
the variable $\delta t$ than all the terms from the former.

Accordingly, in the formally exact result following from \eqref{iterative rule for mu}
all the contributions involving `unknown' things have to be combined with other ones of the same of higher degrees to a single term
to be declared the new
`unknown' contribution devoid of more detailed
internal structure.
It should be of a higher degree in $\delta t$ than all `known' terms ---
otherwise, the accuracy of the input data used in computation is to be declared insufficient 
assuming that, before proceeding, it should be
amended.

In particular, this is the reason why exactly four
transformations are included in the sequence \eqref{eq::720}.
Indeed, the next (the fifth) one,
namely, the transition
$ \{\mu_5,\lamBda_1\} \mbox{\:{\scriptsize$\models$}\:} \mu_6  $,
fails due to insufficient accuracy of the approximation $ \lamBda_1 $, the only one currently available.
As to the transformations \eqref{eq::720}, `the known parts' of their outcome comprise the terms bound with degrees of $\delta t$ up to
1,2,3, and 4,
respectively. The terms of equal degrees from distinct order approximations coincide.
This is the property of the transformations used for their derivation.

To pursue further,
one needs to obtain more accurate (and, accordingly, marked with greater indices) approximations
$ \lamBda_n $. This is realized by means of the formula \eqref{iterative rule for lamIII}.
It involves the kernel function $\Omega_\lamBda$ defined by \Eq{}\eqref{Omega_lam definition}.
In the expression of $\Omega_\lamBda$
$\mu$-function is situated on the topmost level
(in the first line of \eqref{Omega_lam definition}),
being not affected by a small factor.
Accordingly,
it has to be specified with accuracy not smaller than the planned accuracy of the result.

On the other hand,
the leading contribution of
$\lamBda$-function is combined with the factor of $\delta t$.
Hence one unit lower degree of accuracy is here enough (that is also necessary for the making the computation meaningful).

Our current fund of $\mu$-functions  contains $\mu_2, \mu_3, \mu_4, \mu_5$ that
allows us to produce $\lamBda$-functions with the same range of indices.
This   is realized by means of the transformations
\begin{equation}                                    \label{eq::730}
\{\lamBda_1,\mu_2\} \mbox{\,{\scriptsize$|\!\!\models$}\:}  \lamBda_2,
\{\lamBda_2,\mu_3\} \mbox{\,{\scriptsize$|\!\!\models$}\:}  \lamBda_3,
 \{\lamBda_3,\mu_4\} \mbox{\,{\scriptsize$|\!\!\models$}\:}  \lamBda_4,
\{\lamBda_4,\mu_5\} \mbox{\,{\scriptsize$|\!\!\models$}\:}  \lamBda_5,
\end{equation}
where 
 `$\mbox{{\scriptsize$|\!\!\models$}}$'
denotes
the map which applies
the transformation 
\eqref{iterative rule for lamIII}
and then re-organizes its result
isolating higher order `unknown' contributions
in the way outlined above in the description of the map  `{\scriptsize$\models$}'\!.
`The known parts' of  the functions $\lamBda_2, \lamBda_3, \lamBda_4, \lamBda_5 $
thus obtained
contain all the terms
bound with powers of $\delta t$ of the degrees
up to 1,2,3, and 4, respectively.

Now we are in position essentially equivalent to the one that had taken place at the beginning of the computation.
There we had had the $\mu,\lamBda$-functions with values of indices up to $N=1$. 
We managed further to
construct such kind functions with greater indices lifting their maximal value  to $N=5$. 

Now we may assign the role of the original function collection  $\{\mu_1, \lamBda_1  \}$
to  the set  $\{\mu_1,\dots,\mu_N, \lamBda_1\dots,\lamBda_N\}, \; N=5$,
and
carry out the new computations distinct in notations from the above ones only in values of indices labeling
$\mu,\lamBda$-functions
which, as compared to above, will be incremented by 4.

Eight in total new $\mu$-  and $\lamBda$-functions   will then result.
Designating their pair of the maximal order to play the
role of the input data, the application of the eight-step algorithm outlined above can be repeated.

It is evident that there is no place for arising of 
an obstruction
for a subsequent iterative repetition of the above procedure --- besides the finiteness of available computational resources or the like.

Let us notice again that the iterative computations outlined above  do not actually implements the
transformations \eqref{iterative rule for mu}, \eqref{iterative rule for lamIII} rigorously understood.
At each step here not all contributions which they produce
but only ones of a single next degree in $\delta t$ are kept as new and significant.
The higher order terms prove to be summed up with some constituents
previously declared `unknown' and, accordingly,
the result they are incorporated in is also to be regarded as an `unknown' quantity to be finally eliminated from consideration.

The scheme of computation of the right hand sides of formulas \eqref{iterative rule for mu}, \eqref{iterative rule for lamIII}
involving, besides the linear operator of integration, only multiplications and summations, ensures consistency of such a mechanism.
In particular, at each step (i) a new contribution of the order greater by one unit (as compared to the maximum of the available ones)
in $\delta t$ is generated and (ii) all the terms obtained in this way in the preceding iterations remain unchanged.
The computation reproduces in fact the corresponding initial parts of the power series
which arise as the result of application of formulas \eqref{iterative rule for mu}, \eqref{iterative rule for lamIII} carried out
in a rigorous way without any extraneous, even well-founded, modification.

To tentatively characterize
capability 
of the above computational scheme, we display below
the result of its application  for the maximal index of the accuracy equal to 5
(greater accuracies being also feasible but fairly hard for a readable printed presentation).
It is of major interest to inspect
a $\lamBda$-function which is now named $\lamBda_6$. It was obtained that

\vfill
\begin{equation}                                 \label{lam_6 definition}          
\hphantom{.}\hspace{-5.0em}
\begin{aligned}[t]
\lamBda_6(t)=
 & \;
 \lamIII
-
   \frac{\delta t}{4 t_0^2}
      \Big(\chi_\infty + (\sgn\, \chi_0 + 2)\, t_0 \,\lamIII \Big)
\\&
+
   \frac{{\delta t}^2}{20 t_0^2}
     \, \sgn\, \Big(2 + 3 \chi_\infty\, (\chi_0 + \sgn)/t_0 + (5 \chi_0 + 7 \sgn) \lamIII + 6 t_0^2 \,\lamIIIp^2  \Big)
\\&
-
   \frac{{\delta t}^3}{360 t_0^3}
     \Big(
        \chi_\infty\, \big((\chi_0 + \sgn) (9 \chi_0 + 46 \sgn)/t_0 + 90 \sgn\, t_0 \,\lamIII \big)
\\[-1.7ex]& \hspace{4em}
      + 2 (18 \chi_0 + 7 \sgn) + 9 \sgn\, (9 \chi_0 + 11 \sgn) \lamIII + 18 (\chi_0 + 9 \sgn) t_0^2 \,\lamIIIp^2
    \Big)
\\[-1ex]&
+
  \frac{{\delta t}^4}{2520 t_0^5}
  \Big(
   90 \sgn\,t_0\,\chi_\infty^2
\\[-1.1ex]&\hspace{4em}
 +\chi_\infty\,\big( (\chi_0 + \sgn) (91 \chi_0 + 284 \sgn) + 18 (18 \chi_0 + 53 \sgn)\,t_0^2 \,\lamIII\big)
\\&\hspace{4em}
 + 2 \big(97 \chi_0 + \sgn\,(45 \chi_0^2 + 53) \big)\, t_0
\\&\hspace{4em}
 + 36 (11 t_0^2 + 14 \sgn\,\chi_0 + 16) t_0\,\,\lamIII
\\[-1.1ex]& \hspace{4em}
 + 18 (14 \chi_0 + 73 \sgn)\,t_0^3\,\,\lamIIIp^2
 + 108 t_0^5\,\,\lamIIIp^3
\Big)
\\[-1ex]&
-
  \frac{{\delta t}^5}{20160 t_0^6}
  \Big(
 18 (33 \chi_0 + 65 \sgn)\,t_0\,\chi_\infty^2
\\[-1ex]
&\hspace{5em}
  + \chi_\infty\,\big((\chi_0 + \sgn) (830 \chi_0 + 2047 \sgn) + 756 t_0^2
\\
&\hspace{8em}
            + 36 \sgn\,(9 \chi_0^2  + 140 \sgn\,\chi_0 + 257)t_0^2\,\lamIII + 2268 t_0^4\,\lamIIIp^2\big)
\\
&\hspace{5em}
  + 2 (45 \chi_0^3   + 423 \sgn\,\chi_0^2 + 761 \chi_0 + 388 \sgn)\,t_0
\\
&\hspace{5em}
  + 36 \big(100 \sgn\,\chi_0 + 110 + 27 (3 \sgn\,\chi_0 + 4)\,t_0^2\big)\,t_0\,\lamIII
\\
&\hspace{5em}
  + 18 (157 \chi_0 + 620 \sgn)\,t_0^3\,\lamIIIp^2
\\
&\hspace{5em}
  + 108 (\sgn\,\chi_0  + 20)\,t_0^5\,\lamIIIp^3
\Big)
 ,
\\
&\hspace{-2.5em}
\mbox{where the  replacement }
\delta t
 \leftleftharpoons
t-t_0\mbox{ has to be carried out}.
\end{aligned}\hspace{-5em}
\end{equation}

\vfill

\vphantom{.}

\noindent
It could be mentioned that the shown result proves to be less bulky than one might anticipate looking at formulas used for its derivation.


$\lamBda$-function is however  not the final result.
Taking into account \Eq{}\eqref{eq::120}, we can introduce
the function
to be
proposed
as
  the
approximate solution, vanishing at $t=t_0$, to
the \PainIII{} equation \eqref{equation PIII'}
accurate up to terms proportional to $ \delta t^8$, inclusively,
by the following definition 
\begin{equation}                 \label{PIII approx}               
\lambda(t)
\approx
(t-t_0)\sgn
+\half(t-t_0)^2(\sgn-\chi_0)/t_0
+(t-t_0)^3 \,\lamBda_6(t).
\end{equation}

As an independent verification, 
having substituted the above expression into \Eq{}\eqref{equation PIII'},
it was found 
by means of a routine simplification%
\footnote{One should be warned that such a computation seems unlikely to be feasible without
application of the computer algebra tools.}
 that the residual discrepancy
$\ddot\lambda-\mbox{{\sf RHS}[Eq.\eqref{equation PIII'}]} $
in fulfillment of Eq.\eqref{equation PIII'}
proves to be here proportional to $\delta t^7$ --- one unit in the degree is lost due to presence of the vanishing $\lambda$
in denominators of fractions from the equation right hand side.

Removing one or more higher order terms from the expression \eqref{lam_6 definition}, less accurate approximate solutions arise
which verify \Eq{}\eqref{equation PIII'} up to residual discrepancies going to zero as $t\to t_0$
with the corresponding rates of lower degrees.

\subsection{Numerical example}

In addition to the above analytic 
substantiation,
it seems worthwhile to demonstrate a particular lucid
example
confirming relevance
of the formula \eqref{lam_6 definition} and thereby illustrating capability of the method of its derivation.
A collation of a numerical `practically exact' solution to \Eq{}\eqref{equation PIII'}
with results produced by the mentioned approximate analytic solution seems to be a reasonable form of such an illustration.

\begin{figure}
\scalebox{1.22}{
\includegraphics{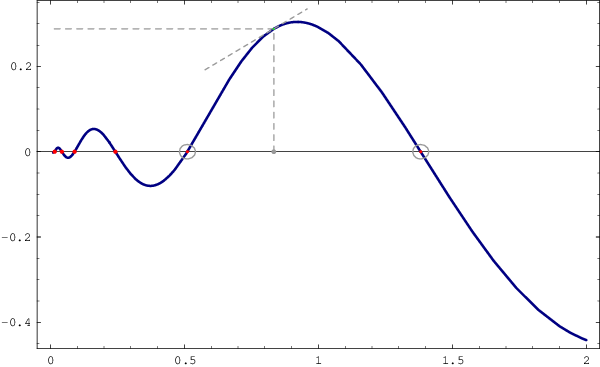}
}
\caption{
The result of numerical integrating of the \PainIII{} equation \eqref{equation PIII'}
on the interval $[0.01,2]$ is shown.
The values of the approximate solution $\lambda(t)$
and of its argument $t$
 are plotted against the vertical and horizontal directions, respectively, 
The constant parameters fixing the equation were set as follows:
$
\chi_0 = -0.811598...$, $\chi_\infty=-0.055004...
$
The solution is fixed
by the values of $\lambda(t)=0.288298...$ and $\dot\lambda(t)=0.374531...$ taken for $t=0.833651... $
that is indicated by the horizontal, inclined, and vertical dashed line segments, respectively,
intersecting at the point of specification of the listed initial data.
}\label{Fig1}
\end{figure}

Dealing with a numerical presentation
a full concreteness is mandatory.
 To that end, we arrange all the constant parameters, free variable, and   functions we will use
to be real valued. We set the constants fixing the equation \eqref{equation PIII'} as follows:
$$
\chi_0 = -0.811597... 
, \chi_\infty=-0.0550042... 
$$
Next we have to construct its solution. This is carried out by means of the numerical solving
the Cauchy problem for some instance of generic initial data. 
We select (originally random) point
$t=t_C$ of its specification and the values (also originally random and then adjusted
by means of a number of unsuccessful attempts) the function to be found and its derivative have to be equal to thereat.
A satisfactory package
of such numerical data
 is as follows
$$
          t_C=0.833651... 
,\; \lambda(t_C)=0.288298... 
,\; \dot\lambda(t_C)=0.374531... 
$$
\Eq{}\eqref{equation PIII'}  is then integrated numerically on the interval $(0.01,2)$.
The result is displayed in figure~\ref{Fig1}.
%
%
%

\begin{figure}
\scalebox{1.22}{
\includegraphics{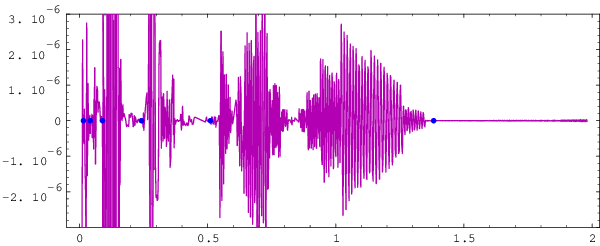}
}
\caption{
The difference of the left and right hand sides of \Eq\eqref{equation PIII'} (the residual discrepancy) 
evaluated on the function $\lambda(t)$ displayed in figure~\ref{Fig1}
is shown. 
The values of the argument
$t$ are plotted against the horizontal coordinate axis. 
The positions of six roots of $\lambda(t)$ are indicated by small disks situated on the latter.
}\label{Fig2}
\end{figure}

The quality of the computed approximation can be indirectly overseen evaluating
the residual discrepancy
$ \ddot\lambda-\mbox{{\sf RHS}[Eq.\eqref{equation PIII'}]}  $,
where
the derivative $ \ddot\lambda = \ddot\lambda(t)$ is determined numerically using an appropriate second order finite difference scheme.
The corresponding graph
exhibiting an apparently chaotic noise caused by
uncontrollable small deviations of the numerical solution from the exact one,
as well as other numerical effects,
is shown in figure~\ref{Fig2}.
It seems indicating, at least for $t>0.3$, a
satisfactory accuracy level.%
%
%

Six zeros are observed on the plot shown in figure~\ref{Fig1}.
Their positions are found numerically as roots of the available approximation $ \lambda(t) $.
The values thus obtained are as follows:
0.0159082...
, 0.0427774...
, 0.0901638...
, 0.242530...
, 0.511115...
, 1.38175... 
On this six-node
grid the approximation  we inspect
takes the following values:
$
-1.7 \cdot10^{-18}
,$ $ 3.5 \cdot10^{-18}
,$ $ 6.3 \cdot10^{-19}
,$ $ 6.5 \cdot10^{-18}
,$ $ 9.1 \cdot10^{-17}
,$ $ 3.3 \cdot10^{-17}.
$
We see that 
the root approximations may be consider satisfactory --- of course, up
to
accuracy of the numerical solution $\lambda$ in itself.
Moreover, as it should be expected,
 at all zeros of $\lambda$
the equality $|\dot\lambda|=1$ takes place with high precision.
Then the
inaccuracy of the above determination of roots has to be the same as
the inaccuracy of the approximate recovering of $\lambda$.

We will more closely consider here only the two maximal zeros from the above their list. 
 In figure~\ref{Fig1}, they are encircled by small circles.

 Near the left root from the pair, $t_{0_1}=0.511115...$, the plot slope is ascending with the derivative $\dot\lambda$
practically coinciding thereat with $+1$.

As it has to be expected,
near the right one, $t_{0_2}=1.38175...$, the plot slope is descending, and at the root the derivative can be get equal to $-1$.
Thus both allowable opportunities $\sgn=\pm1  $ are realized.

We intend to examine capability of the approximated solution \eqref{lam_6 definition} in cases when it is `anchored'
at each of these roots of $\lambda$, independently. 

However, there is an additional problem which must be addressed beforehand. The point is that computations with formula
\eqref{lam_6 definition} require knowledge of one more parameter characterizing solution
to be approximated
 but also depending on the root under consideration.
Above, it was referred to using the notation $\lamIII$.

From a general point of view $\lamIII$ is related to what a concrete solution in their family is dealt with. On the other hand,
in a narrow sense, it
is defined by
the third order derivative of the solution determined at the root at which it is anchored.
Thus, originally, determination of $\lamIII$ is a numerical computation.
At the same time,
the problem
can be lifted  to a partially analytic level if we take for granted that the function to be thrice differentiated
verifies  a second order differential equation. Then one can obtain the third derivative by means of the
differentiating of this equation, i.e.\ differentiating \Eq{}\eqref{equation PIII'}.
Using it once more, the second derivative is eliminated from the
result and one obtains an explicit formula representing the third derivative as a rational function of
the solution itself and its first order derivative (it seems not actually useful to display this formula here).
Both these functions have become known as the result of numerical integration of the original second order equation.

It is not possible, however, to obtain the desirable value of the third order derivative
merely substituting the value of the root into its exact formula noted above since the latter
inherits from \Eq{}\eqref{equation PIII'}
the $0/0$-kind indeterminate behavior
just at values of $t$ where   $\lambda$ vanishes.

A way to surmount that obstacle
could be the using,
instead of
straightforward application of
the derived exact formula,  the
numerical function obtained from the latter by means of interpolation.
Such an interpolation
can be carried out since the third derivative is well defined and smooth (analytic) at roots of solutions, as well
as the lower order derivatives and the solution itself.

\begin{figure}
\scalebox{1.22}{
\includegraphics{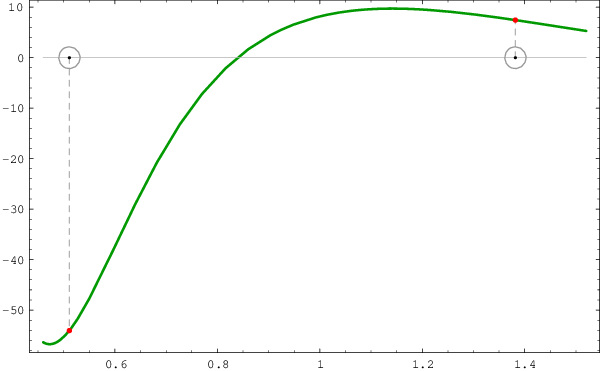}
}
\caption{The third order derivative of the function $\lambda(t)$ displayed in figure~\ref{Fig1} is plotted on a part of its domain.
The positions of the two selected roots of $\lambda(t)$ picked out therein are indicated.}\label{Fig3}
\end{figure}

In practice, the numerical solution we have at our disposal%
\footnote{It was obtained by application of the routine {\sf NDSolve} included in the software {\sc mathematica}.}
is actually represented by a set of its  values and values of its first order derivative at nodes of some finite inhomogeneous grid.
For any point distinct from all the grid nodes (but falling inside the segment the latter span)
the determination of the solution value is realized by means of interpolation (by defaults, using splines).
Working with the third order derivative,
we may proceed in the same way: first of all, its values are determined  at nodes of the noted grid by means of the exact formula
using the stored values of the solution and its derivative at the same node.
It could prove reasonable to remove
some `bad nodes' which turn out to be too close to some of $\lambda$ roots
but in the case we consider no such a necessity arises.
These values are stored.
After that, the third order derivative
is determined at an arbitrary interior point which, in particular,
may coincide with a solution root, by means of interpolation between the nearby grid nodes.
The result of application of such an technique is displayed in figure~\ref{Fig3}. 

\begin{figure}
\scalebox{1.22}{
\includegraphics{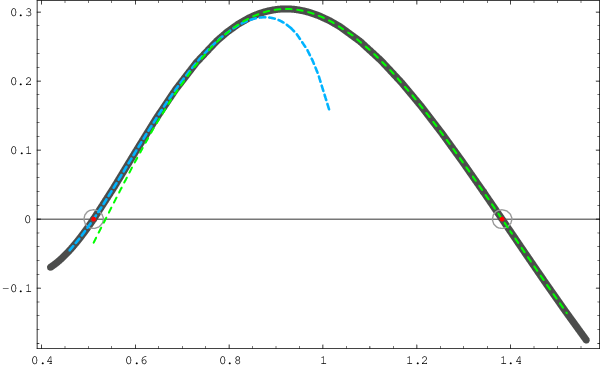}
}
\caption{
Overlapping between
the
two graphs imaging the
approximations defined by the formula \eqref{lam_6 definition}
which are centered at 
$t_{0_1}=0.5111...$ and $t_{0_2}=1.3818... $, respectively, (dashed lines)
and
the plot of numerical solution to \Eq{}\eqref{equation PIII'}
vanishing at $ t_{0_1} $ and $ t_{0_2} $
 (gray solid line in background)
also displayed on a greater domain in figure~\ref{Fig1}.
}\label{Fig4}
\end{figure}

Now,
using the numerical function providing value of the third order derivative of the approximate solution displayed in figure~\ref{Fig1}, we can  determine the values of
$\lamIII_{1}
=\dddot\lambda(  t_{0_1})/6 \approx -9.01149
$  and
$\lamIII_{2}
=\dddot\lambda(  t_{0_2})/6 \approx 1.24246
$ at the roots of $\lambda(t)$ we consider.
As a result,
we have already
 all the data necessary for computations using the formula  \eqref{lam_6 definition} at our disposal.
The values of $\lambda(t)$ obtained by means of its approximate reconstruction on the base of an explicit expression
is displayed in figure~\ref{Fig4}, 
where the interval of variation of $t$
is chosen to be a little wider than the interval $( t_{0_1}, t_{0_2}) $ bounded by the nearby roots selected.

One can see
that the plot of the approximation of the \PIIItran{} we consider by the truncated power series of the order 8 centered at (`anchored' to) the root $t_{0_1}=0.5111...$
looks 
indistinguishable from the plot of the numerical solution
up to the right boundary close to $0.85$ (that corresponds to value of the expansion variable $\delta t/ t_{0_1} \simeq 0.66$).
The left boundary of similar coincidence for the equivalent approximation by the truncated series centered at another, greater, root $t_{0_2}=1.3818...$
is about $0.7$ (corresponding to $\delta t/ t_{0_2}\simeq -0.49$).
Therefore it is situated to the left of the above right boundary of the domain of admissibility of application of the preceding approximation.
Inbetween, around $t\approx 0.8 $, all the three curves depicted in figure~\ref{Fig4}
are most close in total.
Accordingly,
the combination of the above two  approximations
applied on the domains
ensuring their satisfactory accuracies
covers the whole interval of argument variation
bounded by the roots $ t_{0_1} $ and $ t_{0_2} $
and also extends somewhat beyond.

The considered example can be regarded as a numerical validation of the analytical method discussed above and the approximate
formula it has led to.

\section{Application to poles}\label{sec::070}
As it has been mentioned in Introduction, on the set of \PainIII{} transcendents the properties of roots
and poles are closely related.
The details of such a correspondence can be inferred from the existence of the
involutive map
realized by the replacements \eqref{root2pole and vv}
which take solutions to \Eq\eqref{equation PIII'}
to solution of the equation from the same family distinct in values of the constant parameters $\chi_0$ and $\chi_\infty$ which
for the transformed solution
are to be interchanged.
The form of the noted transformation obviously says, in particular, that any non-zero root of the first solution
is situated at the location of a simple pole of the second one, and vice versa.

In view of the above relationships and
taking into  account the assertions of the
theorem  \ref{teo::020},
one can put forward the following statement.
                    \begin{theorem}             \label{teo:160}
Given arbitrary point $t_0\not=0$,
there exist two subfamilies of the family of meromorphic solutions to  \Eq{}\eqref{equation PIII'}
each of which possesses simple pole at $t_0$.
These subfamilies are distinguished by the values of
the residue
which is equal to $t_0$ for one of them and to $-t_0$ for another.
Within a single subfamily,
a solution
can be identified
by
the value of the derivative
$$
{{d\over d t}\:
\bigg\lgroup\!
\lambda(t) - {\mathrm{Res}_{t_0}(\lambda)  \over  t-t_0}
\bigg\rgroup}
$$
evaluated at $t_0$
which is well defined and,
depending on solution,
may manifest any value.
\\
Any singularity of a solution to \Eq{}\eqref{equation PIII'} except the one situated at zero, if any,
represents a simple pole obeying the aforementioned restrictions.
                        \end{theorem}
Taking into account theorem  \ref{teo::020},
the proof reduces to several straightforward computations.\hfill$\square$

The invariance of \Eq{}\eqref{equation PIII'} with respect to the replacements \eqref{root2pole and vv}, \eqref{eq::60} enables one to use the
explicit formula
\eqref{lam_6 definition}
for the constructing of
an approximate representation of \PIIItran{} 
in vicinity of its pole.
The only thing that needs to be here done is the computation of the reciprocal to the given truncated power series.
There are neither principal nor technical obstacles to such a transformation. However, its result proves to be rather
bulky. Albeit being feasible, it seems unlikely to be useful to display it here in full. As a compromise, we provide here less accurate
(or, more precisely, two orders less accurate) version of the noted formula. It reads
%
\begin{equation}                                    \label{eq::760}
\hspace{-4em}
\begin{aligned}[m]
\lambda(t_0+\delta t) = &\;
{ \sgn \,t_0
\over
{\delta t}
}
+
  (\sgn  + \chi_\infty )/2
-
 {\delta t}\,
  \big(\sgn \,(1 - \chi_\infty^2)/(4 t_0)  + t_0\,\lamIII \big)
\\ &
+
 {\delta t}^2/(4 t_0)\,
 \Big(\big((\sgn  - \chi_\infty ) (1 - \chi_\infty^2)/(2 t_0)  + \chi_0\big) + (2 - 3 \sgn \,\chi_\infty ) t_0\,\lamIII
  \Big)
\\ &
-
 {\delta t}^3/(10 t_0)
\times\\ &\hphantom{-}
\Big(    \sgn
   + (3 - 2 \sgn \,\chi_\infty )/(2 t_0)\,\chi_0
   + 5 \big(\sgn \,(1 + \chi_\infty^2) - 2 \chi_\infty \big) (1 - \chi_\infty^2)/(8 t_0^2)
\\ &\hspace{1.8em}
   + \big(1 - 5 (3 \sgn -2 \chi_\infty ) \chi_\infty /2\big)  \lamIII
   - 7 \sgn \,t_0^2\,\lamIIIp^2
\Big)
\\
&
+
 {\delta t}^4/(20 t_0^2)
\times\\ &\hphantom{+}
\Big(
      7 \sgn /9
     + 5 (1 - \chi_\infty^2) \big(\sgn \,(1 + 3 \chi_\infty^2) - (3 + \chi_\infty^2) \chi_\infty \big)/(8 t_0^2)
\\ &\hspace{1.8em}
    + (47 + 45 \chi_\infty^2 - 88 \sgn \,\chi_\infty ) \chi_0/(36 t_0)
\\ &\hspace{1.8em}
    -\big((2 - 15 \chi_\infty^2) + 5 \sgn \,(7 + 5 \chi_\infty^2) \chi_\infty /4 + 5 \sgn \,t_0\,\chi_0\big) \lamIII
\\ &\hspace{1.8em}
    -3 (7 \sgn  - 5 \chi_\infty ) t_0^2\,\lamIIIp^2
\Big)+\Or({\delta t}^5).
\end{aligned}\hspace{-3em}
\end{equation}
A straightforward verification confirms the adequacy of the above  approximation. 
Namely, upon substitution of the approximate solution \eqref{eq::760}
into \Eq{}\eqref{equation PIII'}, the resulting residual discrepancy proves to be proportional to $\delta t^3$.
It could  be added that the making use of the more accurate expression, embodying all the information ensured by \Eq{}\eqref{lam_6 definition},
would impose on the residual discrepancy the overall factor of $\delta t^5$. Thus a producing of such an approximation is quite feasible.
More accurate ones could also be obtained.

\section{Summary}\label{sec::080}

As it is usually assumed
for second order ordinary differential equations,
\PainIII{} equation \eqref{equation PIII'} determines,
almost everywhere, 
the second order derivative $\ddot\lambda$
as a definite regular function (here a rational function) of
the unknown function $\lambda$, its first order derivative $\dot\lambda$,  and their argument $t$.
The obvious exceptions are, at first, the center $t=0$ and, at second, zeros of $\lambda$.
It is
the first order derivative $\dot\lambda$
which is actually determined instead by the equation \eqref{equation PIII'} at roots of $\lambda$.
Moreover, only two values of $\dot\lambda(t_0)$ prove to be there possible, and they are fixed to `universal constants'.
Namely, at a root $t_0$ of  a regular solution $\lambda(t)$ to \PainIII{} equation \eqref{equation PIII'}  (\PIIItran)
it holds either
$\dot\lambda(t_0)=1$
or $\dot\lambda(t_0)=-1$.
Such a peculiarity
means that, searching for \PIIItran{s} vanishing at $t_0\not=0$, it is insufficient to single out them
merely as solutions to the Cauchy problem for \Eq\eqref{equation PIII'} with the requirement $\lambda(t_0)=0$ regarded
as a part of the initial data.
One of the reasons of that is that not only two but more solutions
to the equation vanishing at $t_0$  
should exist but, if this is the case, they cannot be identified by the values of $\dot\lambda(t_0)$.

A better approach 
utilises possibility of regarding \PainIII{} equation as an equation of evolution of some Hamiltonian
dynamical system. The Hamiltonian formalism suggests the form of the replacement of \Eq\eqref{equation PIII'} by a system of two non-linear first order
equations (Hamilton's equations, \Eqs\eqref{Ham eq first}, \eqref{Ham eq second}) which prove to be regular everywhere including roots of $\lambda(t)$
(but still except for $t=0$ where the degeneracy of all the equations retains) and which could be named the coupled Riccati equations.
At this point, the existence of solutions
such that the function $\lambda(t)$ is regular at its roots becomes evident.

It is worthwhile noting that there exist actually (at least) two similar but distinct Hamiltonians leading to \Eq\eqref{equation PIII'}.
They can be defined by the common formula \eqref{Hamiltonian} and arise from it for the distinct values $+1$ and $-1$ of the additional auxiliary
constant parameter $\sgn$. The first of these Hamiltonians is well known while another one seems to be mentioned firstly in Ref.~\cite{T}.

The existence of two Hamiltonians can be associated with existence (or, perhaps better, admissibility) of exactly two
values of $\dot\lambda$ at non-zero roots of  $\lambda$.
Specifically, it turns out that a solution to the system of first order equations
inferred from the Hamiltonian arising when $\sgn=1$ is well defined at a root $t_0$ of $\lambda$ if and only if $\dot\lambda(t_0)=1$.
Analogously, a solution to equations inferred from the Hamiltonian with $\sgn=-1$
is regular at a root $\lambda$ if and only if $\dot\lambda=-1$ thereat.

On the contrary, although the opposite `cross-application' of the noted equation systems in vicinity of roots of $\lambda$
also yields some solutions, their second unknown functions $\mu(t)$ (`generalized momenta') reveal singularities centered thereat.

Avoiding appearance of superfluous irregularities, we impose here the correlation of the
instance of the equations
\eqref{Ham eq first},\eqref{Ham eq second}
(and equivalent \Eqs\eqref{ODE for mu with lamIII}, \eqref{ODE for lamIII})
to be used
with `the discrete part' of the initial data set at a root of solution by means of the link \eqref{eq::100}.

Next, as it is known, considering existence of solutions to differential equations with special properties, it is useful to transform
the latter to appropriate equivalent integral equations. In our case such system
(or two ones, provided the cases of $\sgn=1$ and $\sgn=-1$ are considered independently) consists of
two equations one of which is \Eq\eqref{int eq for mu} while another one can be either \Eq\eqref{int eq for lamIII}
or \Eq\eqref{another int eq for lamIII} which are equivalent. Their solutions can be easily converted to solutions to
\Eq\eqref{equation PIII'} vanishing at the given point.

Using integral equations, the algorithm of development of a sequence of finite power series  of growing degrees
with coinciding starting parts
is formed. See subsection \ref{subs::50.30} of section \ref{sec::040}.
They play role of approximate solutions to the noted integral equations and, at the same time, to differential equations
determining finally \PIIItran{s}.

A particular instance of such an approximate solution is obtained  in section \ref{sec::060}.
Its capability is demonstrated
with the help of a numerical example accompanied with relevant graphical material.
Similar explicit approximation to \PIIItran{} estimating it in vicinity of a pole is also given in section \ref{sec::070}.

The sequence of approximate solutions represented by truncated power series
obtained with the help of the mentioned integral equations
converges uniformly on a certain domain defining therefore an analytic
function which converts to \PIIItran{} by means of a simple transformation.
The necessary substantiation is given in section \ref{sec::040}. It is resumed in theorem \ref{teo::120}
which implies, in particular, the following
                        \begin{corollary}           \label{coro::170}
For any given $t_0\not=0$ and any constant parameter $\lamIII $, complex or real valued,
there exists a solution $\mu(t),\lamBda(t)$
to integral equations \eqref{int eq for mu}, \eqref{another int eq for lamIII}
analytic  in some vicinity of $t=t_0$
which verifies also the differential equations
  \eqref{ODE for mu with lamIII},  \eqref{ODE for lamIII}
and take at $t_0$ the following values 
$$
\lamBda(t_0)=\lamIII,\;
\mu(t_0)  =
\half\big(
1
+
  \sgn\,(1 - \chi_0^2)/(2t_0) + 3 \,t_0 \,\lamIII   \big).
$$
Here $\sgn$ is the additional constant parameter which can be equal to either $+1$ or $-1$.\\
Thereby, in accordance with the foregoing and taking into account propositions \ref{coro::090}, \ref{prop::100},
it can be affirmed that
 the function
\begin{equation*}
\lambda(t) = (t-t_0)\sgn           + \half t_0^{-1}(\sgn - \chi_0)(t-t_0)^2 + (t-t_0)^3 \lamBda(t),
\end{equation*}
vanishing at $t_0$ and analytic thereat, verifies \PainIII{} equation \eqref{equation PIII'}.
                      \end{corollary}
One can elicit from here the degree of arbitrariness in specification of such \PIIItran{s}: they are identified by the sign of $\sgn$
in conjunction with
the value of $\lamIII$. 
The direct products
$\mathbb{Z}_2\times \mathbb{C}$ or $\mathbb{Z}_2\times \mathbb{R}$ can thus be used as the corresponding indexing sets.
From a somewhat distinct perspective the above results are also set out in theorem \ref{teo::020}.

As it has been mentioned in proposition \ref{prop::050},
\PIIItran{s} are analytic at their roots.
At the same time, the above corollary identifies in fact all the analytic solutions to the \PainIII{} equation vanishing at the given point.
Combining these relations, we obtain the following
\begin{proposition}                                \label{prop::180}
As a consequence of the \Pain{} property which \PainIII{} transcendents obey
corollary \ref{coro::170}
specifies all of them which vanish at a given point $t_0\not=0$.
\end{proposition}

The recurrence relations linking elements of sequences of polynomials 
referred to in theorem \ref{teo::120}
can be used for estimating the rates
of their variation that further implies existence of a certain bound for the rate of variation
of the solutions
$\mu(t),\lamBda(t)$.
Basing on such a limitation,
a
bound from below
for the distance which the argument has to pass before \PIIItran{}, having started at definite
point from zero, comes back to zero again, is derived. See theorem  \ref{teo::140}.
The corresponding final estimate for the minimum of distance between two roots of a  \PIIItran{}
can be represented in terms of an elementary function.

 A similar estimate exists for the minimum of distance between the given non-zero root of a \PIIItran{}
and any its pole. Indeed, a pole cannot be situated inside the domain
of convergence of the Taylor series defining \PIIItran{} in question in some vicinity of its root.
At the same time, it follows from theorem \ref{teo::140} that the minimum of maxima of sizes of such domains
can be bounded from below by means of analysis of the noted recurrence relations.
Its estimate is given in corollary  \ref{coro::130}.




\vfill

\appendix

\section{Derivation of integral equations \eqref{int eq for mu}, \eqref{int eq for lamIII}
from differential equations \eqref{ODE for mu with lamIII},\eqref{ODE for lamIII}   }\label{app::A}

As a  straightforward computation shows, the following equalities
\begin{eqnarray}
                                      \label{eq::770}
&&\hspace{-3.em}
\begin{aligned}
{\partial \over \partial \sigma }
\mu(\tau)
\equiv &\;
{\delta t \over t_0}
\bigg(
-
\left\lgroup
{\partial \over \partial \sigma }
\,\sigma\,
\big (
  \half(\chi_\infty+\sgn\,\chi_0-1)+\mu(\tau)
 \big)
\right\rgroup
+
  \Omega_{\mu}[\mu,\lamBda](\sigma,\delta t)
\\& \hphantom{ -{\delta t \over t_0}\Big(  }\!\!
+
\Big(
\tau\,\dot\mu(\tau) -  W_{\!\mu}[\mu,\lamBda](\delta t\,\sigma, \tau)
\Big)
\bigg),
\end{aligned}
\\
                                       \label{eq::780}
&& \hspace{-3.81em}
\begin{aligned}
{\partial \over \partial \sigma }
\sigma^3 \lamBda(\tau)
\equiv &
-{\delta t \over t_0}\left\lgroup
{\partial \over \partial \sigma }
\sigma^4 \lamBda(\tau)
\right\rgroup
+ {\sigma^2\over t_0} \,\Omega_{\lamBda}[\mu,\lamBda](\sigma,\delta t)
\\&
+{\delta t\, \sigma^3  \over t_0}
\Big(
\tau\,\dot\lamBda(\tau) -  W_{\lamBda}[\mu,\lamBda](\delta t\,\sigma, \tau)
\Big),
\end{aligned}
\\ &&\hspace{-3.em}
\text{where $\tau$ has to be replaced with $t_0+(t-t_0)\sigma$ prior to computatuins},
\nonumber
\end{eqnarray}
are fulfilled for arbitrary continuously differentiable functions $\mu,\lamBda$.
One
may
determine `$\sigma$-average' of the left and right hand parts of the both above identities, that is, compute
the result of application to them of the integral operator $\int_0^1\! d\,\sigma\times$.
It has however to be noted that in view of assumed fulfillment of
\Eqs\eqref{ODE for mu with lamIII}, \eqref{ODE for lamIII},
the expressions in parentheses in second lines of their records vanish and these terms may be discarded.
We have to integrate therefore only contributions arising from the first lines of the formulas in question.
On the left, we obviously obtain
 $\mu(t)-\mu(t_0)$  and  $\lamBda(t)$  respectively.
The first summands on the right,
containing the operator $ -(\delta t/t_0)(\partial/\partial \sigma)$,
are reduced
in a similar way
yielding certain explicit expressions.
The integrals of terms  proportional to  $ \Omega_{\mu} $ and $ \Omega_{\lamBda} $
are kept unchanged.
Then one can easily find that the result of the above transformation (`$\sigma$-averaging')
of the identities  \eqref{eq::770}, \eqref{eq::780}
is
just
the equations
\eqref{int eq for mu} and \eqref{int eq for lamIII}, respectively,
provided \Eq{}\eqref{smoothness condition} is taken into account.
The latter follows from the same \Eq\eqref{ODE for lamIII} and the condition of continuity of its right hand side at $t=t_0$, see the first summand in definition
\eqref{W_lam definition}.
\hfill
$\square$

\section{Derivation of \Eq(\ref{smoothness condition}) and \Eq (\ref{ODE for mu with lamIII})
from the integral equation (\ref{int eq for mu})}\label{app::B}

At first,  for $t=t_0$ \Eq\eqref{int eq for mu} takes form of
the equality equivalent to \Eq\eqref{smoothness condition}.

At second, the fulfillment of \Eq(\ref{ODE for mu with lamIII})
follows from the identity
%
\begin{eqnarray}
                                                 \label{eq::790}
&&\  
\begin{aligned}
{\partial \over \partial \sigma }\,
\sigma\,
\Big(
\tau\,\dot\mu(\tau) -  W_{\!\mu}[\mu,\lamBda](\sigma\,\delta t, \tau)
\Big)
\equiv&
\\& \hspace{-15.2em}
{\partial \over \partial t }
\Big(
\left\lgroup
{\partial \over \partial \sigma }
\big( \tau\,\mu(\tau) + \sigma\,\delta t\,(\chi_\infty + \sgn\,\chi_0 -1)/2 \big)
\right\rgroup
- \delta t\,\Omega_{\mu}[\mu,\lamBda](\sigma,\delta t)
\Big),
\\& \hspace{-15.2em}
\mbox{in which the replacements }
\tau \leftleftharpoons\, t_0+ (t-t_0)\sigma,\;
\delta t \leftleftharpoons\, t-t_0
\mbox{ have to be carried}
\\[-1ex]& \hspace{-15.2em}
\mbox{out prior to computations}.
\end{aligned}
\end{eqnarray}
The above equality takes place 
for arbitrary sufficiently smooth (in particular, analytic) functions  $\mu(t), \lamBda(t)$.
$\sigma$-average of its left hand side 
coincides with difference
of the left and right hand side of \Eq(\ref{ODE for mu with lamIII}). It remains to show
that $\sigma$-average of the right hand side of \Eq\eqref{eq::790} vanishes. Adjusting its expression, we note first
that $\partial/\partial t$-derivative is here exchangeable with the integral operator $\int_0^1\! d\sigma\times$. Applying the latter to the target of the former,
we obtain, upon obvious transformation of the first term, the difference of the left and right hand sides of \Eq\eqref{int eq for mu}
times $t_0$. It is assumed to vanish implying thus the expectable vanishing of $\sigma$-average of the right hand side of \Eq\eqref{eq::790}.
Then the fulfillment of \Eq(\ref{ODE for mu with lamIII}) follows.\hfill$\square$
\bigskip

\section{Integral representation of potentially singular term from \Eq\eqref{ODE for lamIII} }\label{app::C}

{\marginparsep=-1cm 
Let us consider $\sigma$-averaging
of the following three expressions
%
%
\\[-1.5ex]
\hphantom{.}\hspace{-5em}
\begin{description}
\item[(A)]
~\\[-2.8em]
\begin{equation}\hspace{-1em}
                                \label{eq::800}
\hphantom{-\,}
{{\delta t}
\over
t_0}
\left(
-\!\!
\left\lgroup
{
\partial
\over
\partial \sigma
}
\sigma
\big(
        \half(\chi_\infty+ \sgn\,\chi_0 - 1)+\mu(\tau)
\big)
\right\rgroup
+\Omega_{\mu}[\mu,\lamBda](\sigma,{\delta t} ) \right),
\end{equation}
~\\[-3ex]
\item[(B)]
~\\[-2.8em]
\begin{equation}\hspace{-7em}
                             \label{eq::810}
-{{\delta t}
\over
t_0}
\!
\left\lgroup
{
\partial
\over
\partial \sigma
}
\sigma^4
\lamBda(\tau)
\right\rgroup
+
{\sigma^2
\over
t_0
}
\Omega_{\lamBda}[\mu,\lamBda](\sigma,{\delta t} ),
\mbox{ and}
\end{equation}
~\\[-2.2ex]
\item[(C)]
~\\[-2.9em]
\begin{equation}\hspace{5em}
                                 \label{eq::820}
-
\left\lgroup
{
\partial
\over
\partial \sigma
}
{\sigma^4
\over
t_0}
\big(
(\chi_\infty+\sgn\,\chi_0-1)/4+ 2\mu(\tau) - 3t_0\,\lamBda(\tau)
\big)
\right\rgroup
-
{\sigma^3\over t_0}
  \Omega_{\xi}[\mu,\lamBda](\sigma,{\delta t} ),
\end{equation}
\end{description}
in which the replacements
$\tau \leftleftharpoons t_0 +(t-t_0)\sigma, \: \delta t \leftleftharpoons t-t_0$
have to be carried out prior to computations.
It is easy to see that
in case of fulfillment of \Eqs\eqref{int eq for mu}, \eqref{int eq for lamIII}
the corresponding values prove to be equal to
\begin{description}
\item[(A)]
\hspace{2em}$ \mu(t)-\mu(t_0) $,
\item[(B)]
\hspace{2em}$ \lamBda(t) $, \mbox{ and}
\item[(C)]
\hspace{2em}$ \upxi(t), $
\end{description}
respectively.
}

Indeed,
in case   {(A)}, having integrated the expression \eqref{eq::800},
one obtains the right hand side of \Eq{}\eqref{int eq for mu}  without the first summand independent of $t$ and equal to $\mu(t_0)$.
Subtracting it from the left hand side, the difference $ \mu(t)-\mu(t_0) $ arises.
In case   {(B)},
the integrating of the expression
\eqref{eq::810} produces the right hand side of \Eq{}\eqref{int eq for lamIII} also fulfilled
and hence its left hand side is equal to $ \lamBda(t) $.
In the last case (C), the integral of the expression \eqref{eq::820} reduces obviously to the right hand side
of the definition \eqref{upxi definition} of the function $\upxi(t)$.

The above notes imply that
in case of fulfillment of \Eqs\eqref{int eq for mu},\eqref{int eq for lamIII}
$\sigma$-average of the linear combination
\begin{equation}
                          \label{V introduction}  
\mathcal{V}=
     2\, \mbox{\sf Expr\eqref{eq::800}}
 - 3 t_0\,\mbox{\sf Expr\eqref{eq::810}}
   -(t-t_0)\,\mbox{\sf Expr\eqref{eq::820}}
\end{equation}
of the expressions referred to by their numbers indicated in this record can be represented in explicit form as follows
\begin{equation}
                                     \label{eq::840}
\int_0^1\!\!\mathcal{V} d \sigma  =
  2\mu(t) - 2\mu(t_0)
- 3 t_0\,\lamBda(t)
- (t-t_0)\upxi(t).
\end{equation}

However,
the original representation \eqref{V introduction} is not an optimal form of $ \mathcal{V} $.
Expanding relevant definitions, one can derive by means of straightforward computation its alternative more sound representation which reads
%
%
\begin{equation}
                                     \label{V definition}   
\hspace{-0em}
\begin{aligned}
\mathcal{V}=&\;
-
{t-t_0
\over
4t_0}
\left\lgroup
{
\partial
\over
\partial \sigma
}\sigma
\left(
3( \chi_\infty+\sgn\,\chi_0-1)+(1-\sigma^3)\big((\chi_\infty+\sgn\,\chi_0-1)+8\mu(\tau)\big)
      \right)
\!\right\rgroup
\\ &\;
-3\sigma^2 \left(
\sgn\,(\chi_0^2-1)/(2t_0) - 1 + 2t_0^{-1}\tau\,\mu(\tau)
           \right)
\\&\;
+2(1-\sigma^3)t_0^{-1}(t-t_0)\,\Omega_{\mu}[\mu,\lamBda](\sigma,t - t_0),
\\&
\mbox{\hspace{-1.7em}\rlap{where the replacing $\tau \leftleftharpoons\, t_0 +(t-t_0)\sigma$ \:to be carried out prior to computations}}
\\[-1ex] & \mbox{\hspace{-1.7em}\rlap{is assumed.}}
\end{aligned}\hspace{-3em}
\end{equation}
We will use it in proof of the following statement
                          \begin{lemma}                 \label{lemm::190}
 Let the functions $\mu(t), \lamBda(t)$ be regular and the integral equation
\eqref{int eq for mu} be fulfilled
in some vicinity of $t=t_0\not=0$. Then $\sigma$-average  $ \int_0^1 \mathcal{V}  d \sigma$
of the expression
defined by \Eq\eqref{V definition}
does not depend on $t$.
                          \end{lemma}
\begin{proof}
The lemma assertion is equivalent to the claim of the vanishing of $t$-derivative of the value of the integral in question.
In turn, this is equivalent to the vanishing of $\sigma$-average
of the $t$-derivative of the expression \eqref{V definition}.
The necessary rearrangement
of the latter
is ensured by
the  identity
%
\begin{equation}
                                \label{eq::860}
              \hspace{-3em}
\begin{aligned}[t]
{
\partial
\over
\partial t
}
\mathcal{V}
%
                           \equiv &\;
-
 \bigg\lgroup\!
{
 \partial
 \over
 \partial t
}
2(1-\sigma^3)
{t-t_0  \over t_0 }
\big(
\tau\,\dot\mu(\tau) -  W_{\!\mu}[\mu,\lamBda]( \tau-t_0  , \tau)
\big)
\!\bigg\rgroup
\\
& \:\: 
+
\!\bigg\lgroup\!
{
\partial
\over
\partial \sigma
}
2(1-\sigma^3)\,\sigma\,\dot\mu(\tau)
\!\bigg\rgroup,
\\ &  \hspace{-5.0em}
\mbox{ where the replacing }\tau \leftleftharpoons\, t_0 +(t-t_0)\sigma
\mbox{ has to be carried out prior to computation}.
\end{aligned}\hspace{-4em}
\end{equation}

The second summand on the right has the null
$\sigma$-average
due to the multiplier  $ (1-\sigma^3)\,\sigma $
vanishing at the boundary points of the integration interval.

In the first summand, the difference $ \tau\,\dot\mu -  W_{\!\mu} $
coincides with the difference of the left and right
sizes of \Eq{}\eqref{ODE for mu with lamIII} evaluated at the point $\tau=t_0+(t-t_0)\sigma$
belonging do the admissible domain (a convex neighborhood of $t_0$).
\Eq{}\eqref{ODE for mu with lamIII}
is there fulfilled
due to fulfillment of the integral equation \eqref{int eq for mu}.
See appendix \ref{app::B}.
Thus the first summand yields no contribution as well.

In total, $\sigma$-average of the right hand side of \eqref{eq::860} proves to be equal to zero implying
vanishing of the derivative of $\sigma$-average of $\mathcal{V}$ which is therefore independent of $t$. The  lemma is proven.
\end{proof}

\Eq\eqref{V definition} evaluated at $t=t_0$ yields
$
\raisebox{-0.8ex}{\scalebox{0.95}{\Big\lfloor}}
\raisebox{-1.6ex}{\mbox{\scriptsize ${t\leftleftharpoons t_0}$ }}
\hspace{-2.2em}\mathcal{V}
=-3\sigma^2 \big(\sgn\,(\chi_0^2-1)/(2t_0) -  1 + 2\mu(t_0)\big)
$. $\sigma$-average of this expression is equal to $ -\big(\sgn\,(\chi_0^2-1)/(2t_0) -  1 + 2\mu(t_0)\big)$.
This is just the value independent of $t$ to which the linear combination \eqref{eq::840} is equal for all $t$ sufficiently close to
$t_0$. In the other words, it holds
\begin{equation}
                                         \label{eq::870}
0
   =
   \sgn\,(\chi_0^2-1)/(2t_0) - 1
+ 2 \mu(t) -   3 t_0\,\lamBda(t)   - (t-t_0)\upxi(t)
\end{equation}
`outside' $t_0$ as well.
This equality has the form of \Eq\eqref{xi fraction} that implies coincidence of functions $\upxi(t)$ and $\xi(t)$. 
We have shown therefore that the fraction \eqref{xi fraction} coincides with regular function arising
as the result of
the integral transformation
\eqref{upxi definition}
of the functions $\mu,\lamBda$.
\hfill$\square$

\section{Derivation of \Eq(\ref{ODE for lamIII}) from integral equations \eqref{int eq for mu}, \eqref{int eq for lamIII} }\label{app::D}

Let us consider the identity
\begin{equation}
                            \label{eq::880}
\hspace{-3em}
\begin{aligned}[t]
3\sigma^3
{\partial\over\partial \sigma}
\sigma
\big(
\tau\dot{\lamBda}(\tau) - W_{\lamBda}[\mu,\lamBda]( \tau-t_0 , \tau)
\big)
\equiv  \hspace{1.6em}  &
\\ & \hspace{-18.7em}
-{\partial\over\partial t}
\bigg(
(t-t_0)
\Big(
\!\Big\lgroup\!
{\partial\over\partial \sigma}
\sigma^4
\big(
\upxi(\tau)
+
 t_0^{-1}
\big(
\quart(\chi_\infty + \sgn\,\chi_0 - 1) + 2\mu(\tau) - 3t_0 \,\lamBda(\tau)
\big)
 \big)
\! \Big\rgroup
\\[-1ex]
&  \hspace{-11em}
+t_0^{-1}\sigma^3\,\Omega_{\xi}[\mu,\lamBda](\sigma, t - t_0)
\Big)
%
\\
  &   \hspace{-16em}
- 2\sigma^3\,t_0^{-1}(t-t_0)
\big(
\tau\,\dot\mu(\tau) -W_{\mu}[\mu,\lamBda](\tau-t_0, \tau)
\big)
\end{aligned}\hspace{-3em}\end{equation}\begin{equation}\nonumber\hphantom{.}\hspace{-5em}\begin{aligned}
&  
+
\Big\lgroup
{\partial\over\partial \sigma}\sigma^3
\left(
\sgn\,(\chi_0^2-1)/(2t_0) - 1 + 2\mu(\tau) - 3t_0\,\lamBda(\tau)- (t-t_0)\sigma\,\upxi(\tau)
\right)
\!\!\Big\rgroup
\!\bigg)
\\&  
 \mbox{in which the replacement }\tau \leftleftharpoons\, t_0+ \sigma\,(t-t_0)\mbox{ has to be carried out prior to}
\\[-1ex]&  
 \mbox{computations.}
\end{aligned}\hspace{-3em}
\end{equation}
%
It holds for arbitrary sufficiently smooth functions $\mu,\lamBda,\upxi$.

As it is shown in appendix \ref{app::C},
in
case of fulfillment of
the integral equations
\eqref{int eq for mu} and \eqref{int eq for lamIII}
the function
 $\upxi$
defined by \Eq\eqref{upxi definition}
obeys \Eq{}\eqref{eq::870}.
Collating the latter with 
\Eq\eqref{eq::880}, one finds that then the last summand in its right hand side vanishes point-wise.
Similarly, the fulfillment of \Eq\eqref{int eq for mu} leads to the fulfilment of \Eq\eqref{ODE for mu  with lamIII}.
This is shown in appendix \ref{app::B}.  Since the multiplier in parentheses in the second summand in the right hand side of
\Eq\eqref{eq::880} is the difference of the left and right sides of \Eq\eqref{ODE for mu  with lamIII} evaluated at the value $\tau$
of the variable $t$ used in the latter, 
the noted summand vanishes point-wise as well.
We are left
therefore
with
the version
of the equality \eqref{eq::880}
retaining only its three topmost lines.

Now let us apply to them the differential operator $(\partial/\partial t)(t-t_0)^4 \times $ and,
afterwards,
the operator of $\sigma$-averaging
$\int_0^1\! d\,\sigma \times  $.
It is claimed that 
the identical zero function arises on the right.

To show this,
let us
apply the operator
$ \int_0^1 d\,\sigma\,(\partial/\partial t)(t-t_0)^4 (\partial/\partial t)(t-t_0)\times $
to the expression separated by large parentheses situated in the second and third lines in \Eq{}\eqref{eq::880}
(as it has just been noted,
the next 
two lines are point-wise null).
Obviously, the integrating and $\partial/\partial t$-derivative commute 
and we may use the equivalent operator
$ (\partial/\partial t)(t-t_0)^4 (\partial/\partial t)(t-t_0)\,\int_0^1 d\,\sigma\times $.
For it, the identical zero arises already at the step of integrating.

Indeed, applying it, the derivative $\partial/\partial \sigma $ in the second line
is annihilated and the $\sigma$-averaging reduces to the replacements
$\sigma \leftleftharpoons\, 1
,
\tau \leftleftharpoons\, t
$. The result of the
$\sigma$-averaging of the contribution from the third line is kept unchanged.
The subsequent
comparison of what has been obtained with the equality (definition) \eqref{upxi definition}
shows
that the former coincides with the difference of the left and right sides of the latter. Hence the result arisen here is equal to zero.

Thus we have shown that
under the conditions assumed
the application of the operator
$\int_0^1\! d\,\sigma(\partial/\partial t)(t-t_0)^4 \times$
to the right hand side of the identity \eqref{eq::880} yields the identically zero  function.
The same must take place for the left hand side, i.e.\ it has to hold
\begin{equation}
                                                            \label{eq::890}
\begin{aligned}[t]
&0=\int_0^1 \!\! d\,\sigma\,
{\partial\over\partial t}
(t-t_0)^4
\sigma^3
{\partial\over\partial \sigma}\,
\sigma\,
\big(
\tau\dot{\lamBda}(\tau) - W_{\lamBda}[\mu,\lamBda]( \tau-t_0 , \tau)
\big)
 ,
\\& 
 \mbox{where  }\tau \mbox{ has to be replaced with } t_0+ \sigma\,(t-t_0) \mbox{ prior to computations}.
\end{aligned}
\end{equation}

To clarify meaning of such a relation, we will use the following statement.
                                     \begin{lemma}             \label{lemm::200}
Let the function $f(t)$ be twice differentiable at and near $t=t_0$.
Then the constraint
\begin{equation}
                                                              \label{eq::900}
\begin{aligned}
&\int_0^1 \!\! d\,\sigma\,
{\partial\over\partial t}
(t-t_0)^4
\sigma^3
{\partial\over\partial \sigma}\,
\sigma
\,f(\tau) 
=0
\end{aligned}
\end{equation}
in which the replacing $ \tau \leftleftharpoons\, t_0+ \sigma\,(t-t_0) $ carried out prior to computations is assumed,
is met
near $t=t_0 $
if and only if
$f(t)\equiv0$.
                                    \end{lemma}
     \begin{proof}
\noindent
It is easy to check that
\begin{equation}\nonumber
\begin{aligned}
{\partial\over\partial t} (t-t_0)^4 \sigma^3 {\partial\over\partial \sigma} \sigma \,f(\tau)
= &\;
 {\partial\over\partial \sigma}  (t-t_0)^3   \sigma^4 {\partial\over\partial t} (t-t_0)\,f(\tau),
\mbox{provided }\tau = t_0+(t-t_0)\sigma.
\end{aligned}
\end{equation}
Let us apply $\sigma$-averaging  $\int_0^1\! d\,\sigma\times $ to the above equality.

On the left, one recognizes the left hand side of \Eq{}\eqref{eq::900}. Hence
it reduces to  the identically zero function.
$\sigma$-averaging of the right hand side must be equal to zero as well.

On the other hand, 
on the right,   $\partial/\partial\sigma$-derivative annihilates with the integration
yielding the result arising in explicit form after the replacement $\sigma \leftleftharpoons 1 $
(assuming the subsidiary replacement  $\tau \leftleftharpoons t $)
everywhere in the part of the expression situated to the right.
In other words, 
the above identity and
\Eq{}\eqref{eq::900}
lead to the equality
\begin{equation}\nonumber
\begin{aligned}
0= (t-t_0)^3   {\partial\over\partial t} (t-t_0)\,f(t).
\end{aligned}
\end{equation}
But this is an indication of the constancy of values of the product $ (t-t_0) f(t)$.
Thus it holds $ f(t) = \mathrm{const}/(t-t_0)$.
However,
such $f(t)$ 
is unbounded in vicinity of the point $t_0$.
Hence 
it is
not continuous and not smooth thereat.
This peculiarity would contradict the lemma conditions
--- unless the constant numerator of the above fraction vanishes.
In the other words,
$f(t)$ must be point-wise equal to zero.
This is just what the lemma assertion means.
\end{proof}

Let us now collate
\Eqs\eqref{eq::890} and \eqref{eq::900}.
Obviously, the former is a particular instance of the latter.
It is important that,
as we had seen,
the expression
$ t\, \dot{\lamBda}( t ) - W_{\lamBda}[\mu,\lamBda](  t -t_0 , t)  $,
playing role of the function $f(t)$,
is a regular function of $t$ including the case of $t=t_0$.
Then the above lemma establishes its identical vanishing at and around $t_0$.
In other words, \Eq\eqref{ODE for lamIII} is there  fulfilled.
\hfill $\square$

\section{Derivation of integral equation \eqref{another int eq for lamIII}
from integral equations
\eqref{int eq for mu}, \eqref{int eq for lamIII}
}\label{app::E}

The implication in question is obtainable from the following identity provable by straightforward computation:
\begin{equation}                               \label{eq::910}
\hspace{-4em}
\begin{aligned}  
&(t-t_0)
\bigg\lgroup
{\partial
\over
\partial t
}
\bigg\lgroup
{\partial
\over
\partial \sigma
}
\left( 
\sigma^3
 \Big(
\lamBda \left( \tau\right)
+
\sigma\,
{t-t_0
\over
t_0}
 \big(  \left( \chi_\infty+\sgn\,\chi_0-1\right)/ \left( 4t_0\right) + \lamBda \left( \tau\right)\big)
\Big)
\right) 
\\[-1ex]
& \hspace{5.5em}
-
  {t-t_0\over 3t_0^2} \: \wideparen{\Omega}_{\lamBda}[\mu,\lamBda](\sigma,t - t_0)
\bigg\rgroup\!\bigg\rgroup
\\                    \llap{$\equiv$}
      &\;
{1 \over (t-t_0)^{2}}
\bigg\lgroup
{\partial
\over
\partial t
}(t-t_0)^{3}
\bigg\lgroup
{\partial
\over
\partial \sigma
}
\left( 
\sigma^3
 \Big(
\lamBda \left( \tau\right)
+
\sigma\,
{t-t_0\over t_0}\,
  \lamBda \left( \tau\right)
\Big)
\right) 
\\[-0.5ex]   &  \hspace{10em}
-
{\sigma^2
\over
t_0
}{\Omega}_{\lamBda}[\mu,\lamBda](\sigma,t - t_0)
\bigg\rgroup\!\bigg\rgroup
\\[-1.2ex]       &
+
\bigg\lgroup
{\partial
\over
\partial \sigma}
\Big( 
{\sigma^3 \over  t_0}
\big( \sgn\,   (\chi_0^2-1)/(2t_0) - 1 + 2\mu(\tau) - 3t_0\lamBda(\tau) -  (\tau - t_0)\,\upxi(\tau)  \big)
\\[-1ex]&\hspace{3em}
-
\Ratio{2}{3}
(1-\sigma^3)\,\sigma\,
{t-t_0
\over
t_0^2
}\,{\Omega}_{\mu}[\mu,\lamBda](\sigma,t - t_0)
\Big) 
\bigg\rgroup
\\
&
+
{t-t_0\over t_0}
\bigg(
\Big\lgroup
{\partial
\over
\partial \sigma}
\sigma^4\,
\Big(
\upxi(\tau)
+
t_0^{-1}
\big(\!
\left(\chi_\infty+\sgn\,\chi_0-1\right)\!/4
+ 2\mu(\tau)
- 3t_0\lamBda(\tau)
\big)
\Big)
\!\Big\rgroup
\\[-1ex]
        &\hspace{4.5em}
+
t_0^{-1}
\sigma^3
{\Omega}_{\xi}[\mu,\lamBda](\sigma,t - t_0)
\bigg),
\\[-1ex] & \hspace{-0.2em}
\mbox{where the replacing }\tau \leftleftharpoons\,  t_0+ (t-t_0)\,\sigma\mbox{ has to be carried out prior to}
\\[-1ex] & \hspace{-0.2em}
\mbox{computations} .
\end{aligned}\hspace{-5em}
\end{equation}
More exactly, the result of its $\sigma$-averaging (application of the integral operator $\int_0^1\! d\sigma\times  $)
will be 
used.
Specifically, it can be proven that $\sigma$-average of the right hand side of \eqref{eq::910}
is the identical zero, provided \Eqs\eqref{int eq for mu}, \eqref{int eq for lamIII} are fulfilled. 
We will do this `summand by summand'.
Inter alia, three ones have to be processed of which the second summand is combined, in turn, from two.

First of all, let us consider the summand with `overall' $\partial/\partial\sigma$-derivative (the fifth and sixth lines of \Eq{}\eqref{eq::910}).
The derivative annihilates with $\sigma$-averaging
and the result appears as the  difference of the
expression to the right evaluated  with  $\sigma=1$ and  $\sigma=0$.
Then the (sub-)summand  containing  $ {\Omega}_{\mu}$ disappears due to the factor of $(1-\sigma^2)\sigma$.
The (sub-)summand
situated in the fifth line
converts, up to the overall factor of $t_0^{-1}$,
to the right hand side of the \Eq{}\eqref{eq::870}.
Under the conditions assumed (the fulfillment of \Eq\eqref{int eq for mu}, see appendix \ref{app::C})
the latter is
fulfilled and thus we obtain the identical zero again.
In total,  $\sigma$-average of the second summand of the right hand side of \eqref{eq::910} vanishes.

$\sigma$-averaging of the third summand of the right hand side of \eqref{eq::910} (its last and last but one lines)
commutes with multiplication by $(t-t_0)/t_0$. Acting to the  expression
on the right of this fraction, $\partial/\partial\sigma$-derivative which it involves annihilates with integrating
reducing to the replacement $\sigma \leftleftharpoons\, 1  $ 1to be applied to its target. $\sigma$-average
of the term involving $ {\Omega}_{\xi}$ is kept without modification.
We obtain precisely the difference of the left and right sides of the
equality \eqref{upxi definition} serving definition
of the function $\upxi(t)$. It is thus equal to zero for all values of $t$.

Finally, we have to establish the vanishing of $\sigma$-average of the summand situated
in the third and fourth lines of \eqref{eq::910}. The integrating commutes with the
differential operator $(t-t_0)^{-2}(\partial/\partial t)(t-t_0)^3$. Acting by $\sigma$-averaging to the target of the latter, we obtain the sum of
two terms. In the first of them,
$\sigma$-averaging annihilates
the $\partial/\partial\sigma$-derivative
and the result arises upon the replacement
$\sigma \leftleftharpoons\, 1  $ (whereas the complementary replacement $\sigma \leftleftharpoons\, 0$ yields zero).
The  form of the $\sigma$-average of the second term involving $ {\Omega}_{\lamBda}$ is kept without
modification. The result thus obtained
coincides with the difference of the
left and right sides of the integral equation \eqref{int eq for lamIII} which is
assumed to be fulfilled.  Thus $\sigma$-average of the term in question also vanishes.

It has been shown therefore  that  $\sigma$-averaging of the right hand side
of the equality \eqref{eq::910} yields the identically zero function of $t$.
 $\sigma$-averaging
of the left hand side 
commutes with
the differential operator
$(t-t_0){\partial/\partial t}$  
and we may pass to consideration of  its action to the target of the latter.
In it,  $\partial/\partial\sigma$-derivative is annihilated
and the computation of the term from the first line reduces to the replacements
$\sigma \leftleftharpoons\, 1,\; \tau \leftleftharpoons t  $.
The form of   $\sigma$-average of the term involving  $ \wideparen{\Omega}_{\lamBda} $
is preserved. 
It is easy to see that
the result thus obtained coincides with the difference
of the left and right sides of  \Eq{}\eqref{another int eq for lamIII} without the constant summand $\lamBda(t_0)$.
The vanishing of
$\sigma$-average of the
right hand side of \Eq{}\eqref{eq::910} established above
means independence of this expression of $t$.
The value of this `constant' can be found by means
of straightforward computation with $t=t_0$. It is trivial and yields $ \lamBda(t_0) $. 
The resulting equality
of the above $t$-independent expression and $ \lamBda(t_0) $
 is equivalent to the equation
 \Eq{}\eqref{another int eq for lamIII} which is thus fulfilled.
Having scrutinized the above reasoning, one finds that
the above conclusion is just a consequence of fulfillment of \Eqs\eqref{int eq for mu}, \eqref{int eq for lamIII}.
\hfill$\square$

\section{Derivation of differential equations (\ref{ODE for mu with lamIII}), (\ref{ODE for lamIII})
from the
alternative system
(\ref{int eq for mu}), (\ref{another int eq for lamIII}) of associated integral equations}\label{app::F}

It is useful to emphasize a separate role of \Eq\eqref{int eq for mu} in the relations we discuss.
Namely, in the case $t=t_0$, it reduces to \Eq\eqref{smoothness condition}
binding the values of $\mu(t_0)$ and $\lamBda(t_0)$,  the constituents of `the initial data' for differential equations
we consider.
In turn, the fulfilment of this constraint
implies the regularity 
of the function $W_{\lamBda}[\mu,\lamBda](t-t_0, t)$ involved in \Eq{}\eqref{ODE for lamIII}
which otherwise would be unbounded in vicinity of $t=t_0$.

Thus if \Eq{}\eqref{int eq for mu} is fulfilled then the expressions of
\Eqs\eqref{ODE for mu with lamIII}  and \eqref{ODE for lamIII}
are  regular functions of $t$ regardless of their fulfillment.
Moreover, in accordance with appendix \ref{app::B} 
the fulfillment of \Eq\eqref{int eq for mu} is already sufficient 
for fulfillment of 
\Eq{}\eqref{ODE for mu with lamIII}.
Thus it remains to prove the fulfillment of \Eq{}\eqref{ODE for lamIII}.

\smallskip

To that end, let us consider the identity 
%
\begin{equation}                              \label{eq::920}
\hphantom{.} \hspace{-8em}
\begin{aligned}
\sigma^3
{\partial\over\partial \sigma}\Big(
\sigma\,
\big(
\tau\,\dot{\lamBda}(\tau) - W_{\lamBda}[\mu,\lamBda](\tau -t_0, \tau)
\big)
\Big)
 \equiv&
\\  & \hspace{-17.7em}
{\partial\over\partial t}
\bigg(
\!\Big\lgroup
{\partial\over\partial \sigma}
\Big(
\sigma^3
\big(
          \tau\lamBda(\tau)
        +\sigma\,(t-t_0)(\chi_\infty+\sgn\,\chi_0-1)/(4t_0)
\\[-1.7ex] & \hspace{-13,em}\hphantom{ \sigma^3 \big( }
 -  (\sgn\,(\chi_0^2-1)/(2t_0) - 1)/3
\big)
\\
& \hspace{-16em}
 \hspace{3em}
+\Ratio{2}{3}(1-\sigma^3)\,\mu(\tau)
\Big)
\!\Big\rgroup
\\[-0.ex]& \hspace{-16em}
-
  {t-t_0 \over 3t_0}\,\wideparen{\Omega}_{\lamBda}[\mu,\lamBda](\sigma,t-t_0)
\\ & \hspace{-16em}
-\Ratio{2}{3}(1-\sigma^3)
{t-t_0 \over t_0} 
\big(
\tau\,\dot{\mu}(\tau) - W_{\mu}[\mu,\lamBda](\tau -t_0, \tau)
\big)
\\
& \hspace{-16em}
+
\bigg\lgroup
{\partial\over\partial \sigma}
\Big(
\Ratio{2}{3}(1-\sigma^3)\,\sigma\,
{t-t_0 \over t_0} 
\,
\big(
         (\chi_\infty+\sgn\,\chi_0-1)/2 + \mu(\tau)
\big)\Big)\!\bigg\rgroup  \!
\bigg),
\\ & 
    \hspace{-19.0em}
\mbox{where the replacing }\tau \leftleftharpoons\, t_0 +(t-t_0)\sigma
\mbox{ has to be carried out prior to computation}.
\end{aligned}\hspace{-8em}
\end{equation}

\smallskip

\noindent
It can be shown by straightforward computation that the equality \eqref{eq::920} 
takes place for arbitrary twice differentiable functions $\mu, \lamBda  $.

First of all, we have to
show that
$\sigma$-average
of the right hand side here
vanishes. 
Since this operation (i.e.~the integrating $\int_0^1\! d\sigma\times$)  commutes with
$\partial/\partial t$-derivative  we may apply the latter directly to the expression to which the derivative acts.

It consists of three summands.

The last of them is a derivative itself, this time with respect to $\sigma$.
It is annihilated with $\sigma$-averaging yielding zero due to the  multiplier $(1-\sigma^3)\,\sigma  $
vanishing at boundaries of the integration interval.

The next in order from below is the summand containing the multiplier
$ \big(
\tau\,\dot{\mu}(\tau) - W_{\mu}[\mu,\lamBda](\tau -t_0, \tau)
\big) $ (the fifth line of \eqref{eq::920}).
It
coincides with the difference of the left and right sides of \Eq{}\eqref{ODE for mu with lamIII}
in which the variable $t$ is replaced by the expression $ t_0 +(t-t_0)\sigma$.
It has been noted that this equation is fulfilled at, and in vicinity of, $t=t_0$.
The expression to be averaged is therefore equal to zero point-wise. Its contribution vanishes. 

A contribution
from the summand displayed in the lines number two, three, and four
consists of   two sub-summands.

\vfill \vphantom{.}

The second of them contains the function $ \wideparen{\Omega}_\lambda $.
It is left without modification except for `the ascent' of the factor of $(t-t_0)/(3t_0)$
independent of $\sigma$ through the action of $\sigma$-averaging.

The first sub-summand  is the result of application  of $\partial/\partial \sigma  $-derivative
to the shown explicit expression.
The derivative is annihilated with the the integral of $\sigma$-averaging
yielding 
difference of the results of the replacements
$\sigma   \leftleftharpoons 1   $ and  $\sigma   \leftleftharpoons 0 $ in the $\sigma$-derivative target.
It is equal to
%
%
\begin{equation}                         \label{eq::930}
\begin{aligned}
&          t_0 \,\lamBda(t) - t_0\,\lamBda(t_0) + (t-t_0) \big( (\chi_\infty + \sgn \chi_0 - 1)/(4t_0) +  \lamBda(t) \big)
\\&       -
          \Ratio{1}{3}\big(\sgn\, (\chi_0^2-1)/(2t_0) - 1 + 2 \mu(t_0) - 3 t_0\, \lamBda(t_0)\big).
\end{aligned}
\end{equation}
The second line here vanishes in view of \Eq\eqref{smoothness condition}.
The residual expression, 
when
incorporated with the integral of $\sigma$-averaging
of the term involving $ \wideparen{\Omega}_\lambda $,
proves to be equal, up to the overall factor of $t_0$, to the difference of the left and right sides of the integral equation \eqref{another int eq for lamIII}
assumed to be fulfilled.
Thus the summand we consider also yields no contribution. 

Combining the above conclusions, one sees that  $\sigma$-average of the target of $\partial/\partial t$-derivative in the right hand side of
\Eq\eqref{eq::920} vanishes as well as $\sigma$-average of the very right hand side.


Applying the operator from the right hand side of the mixed operator identity
\begin{equation*}                            
\int^1_0
\!d \sigma
{\partial
\over
\partial t
}(t-t_0)^4\times
\equiv
{\partial
\over
\partial t
}(t-t_0)^4
\int^1_0\!\! d \sigma
\times 
\end{equation*}
to the right hand side of \eqref {eq::920}, one obtains
the operator
$(\partial/\partial t)(t-t_0)^4\times$
acting to
$\sigma$-average of something finally found to be equal to zero. 

Acting by the left hand side operator to the left hand side of \eqref {eq::920},
the left hand side of
\Eq\eqref{eq::900} with $f(\tau)=\tau\,\dot{\lamBda}(\tau) - W_{\lamBda}[\mu,\lamBda](\tau -t_0, \tau)$ arises.
As it had been noted, the function $ W_{\lamBda}[\mu,\lamBda]$ is regular on solutions to \Eq(\ref{int eq for mu}).
Then lemma \ref{lemm::200} from appendix \ref{app::D} can be applied. It yields $f(t)=0$ or, in other words, states that \Eq\eqref{ODE for lamIII} is fulfilled.

Since the fulfillment of \Eq(\ref{ODE for mu with lamIII}) has been established in appendix \ref{app::B} the proof is completed.%
\hfill$\square$

\end{document}